\RequirePackage{fix-cm}
\documentclass[twocolumn]{svjour3}

\smartqed

\usepackage[english]{babel}
\usepackage[T1]{fontenc}
\usepackage[ansinew]{inputenc}
\usepackage{lmodern}
\usepackage{hyperref}
\usepackage{graphicx}
\usepackage{graphics}
\usepackage{xcolor}
\usepackage{colortbl}
\usepackage{subfig}
\usepackage{algorithm}
\usepackage{algorithmic}
\usepackage{amsmath}
\usepackage{amsfonts}
\usepackage{amssymb}
\usepackage{booktabs} 
\usepackage{hhline}
\usepackage{etoolbox}
\usepackage{multirow}
\usepackage{float}
\usepackage{natbib}
\usepackage{tikz} 
\usepackage{hhline}
\usepackage{nomencl}
\makenomenclature

\addto\extrasenglish{
    
}
\addto\extrasenglish{
    
}

\definecolor{LightCyan}{rgb}{0.88,1,1}

\usepackage{listings}
\definecolor{LightCyan}{rgb}{0.88,1,1}
\definecolor{mygreen}{RGB}{28,172,0} 
\definecolor{mylilas}{RGB}{170,55,241}
\newcommand*\circled[1]{\tikz[baseline=(char.base)]
{\node[shape=circle,draw,inner sep=1pt] (char) {#1};}}

\journalname{Journal Name}

\begin{document}

\title{Topology Optimization with linearized buckling criteria in 250 lines of Matlab}

\titlerunning{Topology Optimization with linearized buckling criteria in 250 lines of Matlab}

\author{Federico Ferrari \and
        Ole Sigmund \and
        James K. Guest}

\authorrunning{F. Ferrari \and
			et al.}

\institute{F. Ferrari \and James K. Guest \at
		Department of Civil and Systems Engineering \\
		Johns Hopkins University \\
		Latrobe Hall, 21218 Baltimore MD, USA \\
		\email{fferrar3@jhu.edu; jkguest@jhu.edu}
		\and O. Sigmund \at
		Department of Mechanical Engineering \\
		Technical University of Denmark \\
		Nils Koppels All\'{e} 404, 2800 Kongens Lyngby, Denmark \\
        \email{sigmund@mek.dtu.dk}
        }

\date{Received: date / Accepted: date}

\date{\textbf{Accepted paper, soon to appear in Structural and Multidisciplinary Optimization}}

\maketitle

\begin{abstract}
We present a 250 line Matlab code for topology optimization for linearized buckling criteria. The code is conceived to handle stiffness, volume and Buckling Load Factors (BLFs) either as the objective function or as constraints. We use the Kreisselmeier-Steinhauser aggregation function in order to reduce multiple objectives (viz. constraints) to a single, differentiable one. Then, the problem is sequentially approximated by using MMA-like expansions and an OC-like scheme is tailored to update the variables. The inspection of the stress stiffness matrix leads to a vectorized implementation for its efficient construction and for the sensitivity analysis of the BLFs. This, coupled with the efficiency improvements already presented by \citet{ferrari-sigmund_20b}, cuts all the computational bottlenecks associated with setting up the buckling analysis and allows buckling topology optimization problems of an interesting size to be solved on a laptop. The efficiency and flexibility of the code is demonstrated over a few structural design examples and some ideas are given for possible extensions.
\end{abstract}

\section{Introduction}
 \label{Sec:intro}

We present a compact code aimed at reducing the burden of a computationally intensive task: topology optimization considering buckling performed in Matlab.

Topology Optimization (TO) is actively spreading across several engineering fields, offering innovative solutions to more and more design problems. Most likely, one of the driving forces of this rapid spreading is the availability of much educational software. We believe that educational software is a higly valuable tool, offering a foundation for the understanding of the basics of a research topic and providing a basis for the implementation (and improvement) of methods. Since the forefather 99 line code by \citet{sigmund_01a}, many others educational codes have appeared (see \cite{ferrari-sigmund_20b} for a list), almost all addressing the prototypical problems of compliance (thermal or elastic) design and mechanism design. In recent years, more TO problems have been translated into educational software, addressing different physics \citep{xia-breitkopf_15a}, design parametrizations \citep{yago-etal_20a}, and sometimes making use of external FEM libraries for the most computationally intensive tasks.

Confining ourselves to density-based TO \citep{book:bendsoe-sigmund_2004} and to structural applications, almost all educational codes deal with compliance design and we observe the lack of any software addressing buckling. This is understandable, as buckling TO is still an advanced topic, presenting several difficulties \citep{bruyneel-etal_08a,ferrari-sigmund_19a}, even in its linearized formulation \citep{book:crisfield91}. Within the Matlab enviroment, buckling analysis is computationally demanding not only due to the eigenproblem solution (progress has been made by using multilevel solvers for this \citep{ferrari-sigmund_20a}), but also due to the stress stiffness matrix setup. The latter depends on the local stress distribution, and therefore its customary setup requires potentially time consuming elementwise operations. Based on the direct inspection of the stress stiffness operator, we provide an almost fully vectorized implementation eliminating all the computational bottlenecks associated with the setup of the buckling problem and with the sensitivity analysis of the Buckling Load Factors (BLFs).

The code we present, named \texttt{topBuck250}, solves TO problems involving a combination of compliance, volume and BLFs. Therefore, it stands as a quite comprehensive design tool, accounting for the most relevant structural responses. In particular, the code in \autoref{App:matlabCodeBuckling} is ready to tackle the following problems
\begin{enumerate}
 \item Maximization of the fundamental BLF of a design, subject to compliance and volume constraints,
 \item Minimum volume design, subject to compliance and buckling constraints,
\end{enumerate}
that are highly meaningful for structural designers.

We aggregate multiple objectives or constraints by using the Kreisselmeier-Steinhauser (KS) function \citep{kreisselmeier-steinhauser_79a}, which has proven very robust when dealing with many active constraints \citep{kennedy-hicken_15a,ferrari-sigmund_19a}. The single objective, single constraint optimization problem is then solved by a sequential approximation approach, using monotonic MMA-like approximations \citep{svanberg_87a} and an OC-like scheme to update the variables. For this, we provide a very compact re-design routine, tailored for solving the specific problems here discussed.

The present code builds on the recently published \texttt{top99neo} code \citep{ferrari-sigmund_20b} and inherits all its basic speedups (i.e. fast assembly implementation, use of volume-preserving filters, etc.). Thus, the solution of the eigenproblem will absorb the majority of the computational time and a highly efficient eigensolver \citep{dunning-etal_16a,ferrari-sigmund_20a} could be used to improve this operation. We do not provide here a complete 3D implementation, because \textit{the solution} of the 3D eigenproblem would quickly become prohibitive with the built-in Matlab tools. Nevertheless, all the methods discussed here are easily extendable to 3D and we will show the huge cut in the CPU time they potentially bring in that setup. Also, the methods presented can be extended to some other discretizations, e.g. higher order isoparameteric elements, some mixed elements and level sets.

The paper is organized as follows. In \autoref{Sec:Formulations} we set the stage for density-based TO, introducing the phyisical responses of interest and their sensitivities. \autoref{sSec:SolutionScheme} introduces the optimization problems considered and qualitatively describes the solution scheme. \autoref{Sec:codeStructure} gives an overview of the organization of the \texttt{topBuck250} code and describes the input parameters. The efficient construction of the stress stiffness matrix and the corresponding sensitivity analysis of the BLFs are presented in \autoref{Sec:BucklingFast}. \autoref{sSec:stateAdjointSolvers} gives some information about how the state and adjoint equations are solved. Two examples, immediately replicable with the provided code, are presented in \autoref{Sec:examples} and final remarks and possible extensions are listed in \autoref{Sec:conclusions}. In \autoref{App-notesDiscretizationOperators} we explain how some discretization operators have been compactly implemented in the code and \autoref{App:ksProperties&OCupdate} describes the MMA-like approximation of the local optimization problem and the OC-like update rule. \autoref{App:matlabCodeBuckling} provides the Matlab code.

\section{Setup and formulations}
 \label{Sec:Formulations}

We consider a structured discretization $\Omega_{h}$ of $m$ equi-sized elements $\Omega_{e}$ for a total of $n$ Degrees of Freedom (DOFs). The pseudo-densities $\hat{\mathbf{x}} = \{ \hat{x}_{e} \}^{m}_{e = 1}$, are partitioned between the sets $\hat{\mathbf{x}}_{\mathcal{P}_{0}}$, $\hat{\mathbf{x}}_{\mathcal{P}_{1}}$, describing passive regions where $\hat{x}_{e} = 0$ and $\hat{x}_{e} = 1$, respectively, and the active variables set $\hat{\mathbf{x}}_{\mathcal{A}}$ \citep{ferrari-sigmund_20b}. For the latter, $\hat{x}_{e}$ is linked to the element-based design variables $\mathbf{x}=\{x_{e}\}^{m}_{e=1}$ by \citep{wang-etal_11a}
\begin{equation}
 \label{eq:projectionOperator}
  \hat{x}_{e} = \frac{\tanh(\beta\eta) + \tanh(\beta(\tilde{x}_{e}-\eta))}{\tanh(\beta\eta)+\tanh(\beta(1-\eta))}
\end{equation}
where $\eta\in[0,1]$, $\beta\in[1,\infty)$ and $\tilde{x}_{e}$ is obtained through the linear density filter with minimum radius $r_{\rm min} > 0$ \citep{bourdin_01a,bruns-tortorelli_01a}
\begin{equation}
 \label{eq:denistyFilter}
  \tilde{x}_{e} =
  \frac{\sum^{m}_{i=1}h_{e,i}x_{i}}
  {\sum^{m}_{i=1}h_{e,i}}
\end{equation}
where $h_{e,i}=\max(0, r_{\rm min} - {\rm dist}(\Omega_{i}, \Omega_{e}))$, forall $i, e \in [1, m]$.

Let the Young's modulus be parametrized by the SIMP interpolations with penalizations $p_{K}$, $p_{G}>1$
\begin{equation}
 \label{eq:stiffnessAndStressInterpolations}
  \begin{aligned}
  E_{K}( \hat{x}_{e} ) & = E_{\rm min} + (E_{0} - E_{\rm min}) \hat{x}^{p_{K}}_{e} \\
  \quad
  E_{G}( \hat{x}_{e} ) & = E_{0} \hat{x}^{p_{G}}_{e}
  \end{aligned}
\end{equation}
for stiffness ($E_{K}$) and stress ($E_{G}$), respectively. As usual, $E_{0}$ is the Young's modulus of the solid material and $E_{\rm min} \ll E_{0}$ that of the void.

We consider the following physical quantities: total volume fraction ($|\Omega_{e}|=\mathrm{constant}$)
\begin{equation}
 \label{eq:StructuralVolume}
  f\left( \hat{\mathbf{x}} \right) =
  \frac{1}{m}
  \sum^{m}_{e=1} \hat{x}_{e}
  \: ,
\end{equation}
linearized compliance
\begin{equation}
 \label{eq:StaticCompliance}
  c(\hat{\mathbf{x}}) =
  \mathbf{F}^{T}\mathbf{u}\left(
   \hat{\mathbf{x}}\right)
   \: ,
\end{equation}
and fundamental Buckling Load Factor (BLF)
\begin{equation}
 \label{eq:minBLFcharRayleigh}
  \lambda_{1}\left( \hat{\mathbf{x}},
  \mathbf{u} \right) := 
  \min\limits_{\mathbf{v}\in\mathbb{R}^{n},
  \mathbf{v}\neq\mathbf{0}}
  - \frac{\mathbf{v}^{T}K\left[ 
  \hat{\mathbf{x}} \right]\mathbf{v}}
  {\mathbf{v}^{T}G\left[ \hat{\mathbf{x}},
  \mathbf{u} \right]\mathbf{v}}
\end{equation}

The latter two both refer to the same load vector $\mathbf{F}\in\mathbb{R}^{n}$, that is assumed to have \textit{fixed} direction; therefore, we only consider positive BLFs. $G\left[ \hat{\mathbf{x}}, \mathbf{u} \right]$ is the stress stiffness matrix, depending on the displacement field $\mathbf{u}(\hat{\mathbf{x}}) = K[\hat{\mathbf{x}}]^{-1}\mathbf{F}$, and $K[\hat{\mathbf{x}}]$ is the linear stiffness matrix.

In the following we will actually compute the quantities $\mu_{i} = 1/\lambda_{i}$, which are the eigenvalues of
\begin{equation}
 \label{eq:eigenvalueEquation}
  \left( G\left[ \hat{\mathbf{x}}, \mathbf{u} \right] + \mu K\left[ \hat{\mathbf{x}} \right] \right) \boldsymbol{\varphi} = \mathbf{0}
  \: , \qquad \boldsymbol{\varphi} \neq \mathbf{0}
\end{equation}
such that $\mu_{1}=\max\mu_{i}=1/\min\lambda_{i} = 1/\lambda_{1}$. Then, the max operator is approximated by the smooth aggregation function \citep{kreisselmeier-steinhauser_79a}
\begin{equation}
 \label{eq:ksAggregationLambda1}
   J^{KS}\left[ \mu_{i} \right](\hat{\mathbf{x}}) = 
   \mu_{1}(\hat{\mathbf{x}}) + 
   \frac{1}{\rho} \ln\left( \sum^{q}_{i=1}
   e^{\rho\left( \mu_{i}(\hat{\mathbf{x}}) - 
   \mu_{1}(\hat{\mathbf{x}}) \right)} \right)
\end{equation}
depending on $\rho \in [1, \infty)$ and giving an upper bound on $\mu_{1}$ (thus, a lower bound on $\lambda_{1}$).

The derivatives of the compliance and volume fraction with respect to the pseudo-densities $\hat{x}_{e}$ are
\begin{equation}
 \label{eq:sensitivityComplianceVolume}
  \frac{\partial c}{\partial \hat{x}_{e}} =
  - \mathbf{u}^{T}\frac{\partial K}{\partial\hat{x}_{e}}\mathbf{u}
  \delta_{e\mathcal{A}} \: , \qquad
  \frac{\partial f}{\partial \hat{x}_{e}} = m^{-1}
  \delta_{e\mathcal{A}}
\end{equation}
where $\delta_{e\mathcal{A}} = 1$ if $e\in\mathcal{A}$ and $0$ otherwise. The derivatives of the $i$-th eigenvalue reads \citep{rodrigues-etal_95a}
\begin{equation}
 \label{eq:sensitivityMu}
 \begin{aligned}
  \frac{\partial \mu_{i}}{\partial \hat{x}_{e}} = 
  -[
  \underbrace{\boldsymbol{\varphi}^{T}_{i}
  \frac{\partial G}{\partial \hat{x}_{e}}
  \boldsymbol{\varphi}_{i}}_{\circled{1}} + \mu_{i}
  \underbrace{\boldsymbol{\varphi}^{T}_{i}
  \frac{\partial K}{\partial \hat{x}_{e}}
  \boldsymbol{\varphi}_{i}}_{\circled{2}} -
  \underbrace{\mathbf{w}^{T}_{i}\frac{\partial K}{\partial \hat{x}_{e}}
  \mathbf{u}}_{\circled{3}}
  ]
 \end{aligned}
\end{equation}
where $\mathbf{w}_{i} = K^{-1}[\boldsymbol{\varphi}^{T}_{i}(\nabla_{\mathbf{u}}G)\boldsymbol{\varphi}_{i}]$ is the adjoint vector, and the derivative of \eqref{eq:ksAggregationLambda1} is \citep{raspanti-etal_00a}
\begin{equation}
 \label{eq:ksAggregationSens}
  \frac{\partial J^{KS}[\mu_{i}]}{\partial\hat{x}_{e}}
   (\hat{\mathbf{x}}) =
   \frac{\sum^{q}_{i=1} e^{\rho(\mu_{i}(\hat{\mathbf{x}})-\mu_{1}(\hat{\mathbf{x}}))}
   \frac{\partial \mu_{i}}{\partial \hat{x}_{e}}
   (\hat{\mathbf{x}})}
   {\sum^{q}_{i=1} e^{\rho(\mu_{i}(\hat{\mathbf{x}})-\mu_{1}(\hat{\mathbf{x}}))}}
\end{equation}

We remark that \eqref{eq:ksAggregationSens} is smooth even at design points where some of the $\partial_{e}\mu_{i}$ are not \citep{gravesen-etal_11a}, i.e. points where $\mu_{i}$ are repeated \citep{seyranian_94a}.

With \eqref{eq:sensitivityComplianceVolume} and \eqref{eq:ksAggregationSens} in hand, sensitivities with respect to the design variables $x_{e}$, are given by the chain rule
\begin{equation}
 \label{eq:chainRuleSensitivities}
 \frac{\partial (\cdot)}{\partial x_{e}} = 
 \frac{\partial (\cdot)}{\partial \hat{x}_{e}}
 \frac{\partial\hat{x}_{e}}{\partial \tilde{x}_{e}}
 \frac{\partial \tilde{x}_{e}}{\partial x_{e}}
\end{equation}
$(\cdot)$ may represent $c$, $f$ or $J^{KS}$ and the rightmost terms involve derivatives of the projection \eqref{eq:projectionOperator} and filtering \eqref{eq:denistyFilter} operators \citep{guest-etal_04a,sigmund_07a}.

\subsection{Optimization problems and solution approach}
 \label{sSec:SolutionScheme}

The code in \autoref{App:matlabCodeBuckling} solves the following problems: (1) BLF maximization with compliance and volume constraints (recall that $J^{KS}_{0}[\mu_{i}]$ is linked to $1/\lambda_{1}$),
\begin{equation}
 \label{eq:optProblemMaxBLF}
  \begin{cases}
  & \min\limits_{\mathbf{x}\in\left[ 0, 1 \right]^{m}}
  J^{KS}_{0}[\mu_{i}]\left( \hat{\mathbf{x}} \right) \\
  {\rm s.t.} & 
  J^{KS}_{1}[g_{V}, g_{c}]\left( \hat{\mathbf{x}} \right) \leq 0
  \end{cases}
\end{equation}
and (2) volume minimization with compliance and BLF constraints
\begin{equation}
 \label{eq:optProblemMinVolcstrBLFandC}
  \begin{cases}
  & \min\limits_{\mathbf{x}\in\left[ 0, 1 \right]^{m}}
  f\left( \hat{\mathbf{x}} \right) \\
  {\rm s.t.} &
  J^{KS}_{1}[g_{c}, g_{\lambda}]\left( \hat{\mathbf{x}} \right) \leq 0
  \end{cases}
\end{equation}
where we used the KS function to aggregate multiple constraints (see \autoref{App:ksProperties&OCupdate}). Explicitly, the constraint functions that are aggregated in \eqref{eq:optProblemMaxBLF} and \eqref{eq:optProblemMinVolcstrBLFandC} are
\begin{equation}
 \label{eq:constraintsForm}
 \begin{aligned}
  g_{V}(\hat{\mathbf{x}}) & = f(\hat{\mathbf{x}}) / \bar{f} - 1 \\
  g_{c}(\hat{\mathbf{x}}) & = c(\hat{\mathbf{x}}) / \bar{c} - 1 \\
  g_{\lambda}(\hat{\mathbf{x}}) & = 1 - \underline{\lambda} J^{KS}[\mu_{i}](\hat{\mathbf{x}})
 \end{aligned}
\end{equation}
where $\bar{c}$ and $\bar{f}$ are the maximum allowed compliance and volume fraction and $\underline{\lambda}$ is the minimum prescribed BLF. The two more basic TO problems of volume-constrained minimum compliance and \emph{vice versa} can also be solved (see \autoref{Sec:codeStructure}).

The optimization problems are solved by a sequential approximation approach \citep{groenwold-etman_08a}. At each re-design step, the objective and constraint functions are replaced by their monotonic MMA-like approximations \citep{svanberg_87a} and the design update is performed within the routine \texttt{ocUpdate} listed in \autoref{App:ksProperties&OCupdate}. Although introducing moving asymptotes and adaptive move limits results in some extra parameters and operations, for the present problems the design update retains the simplicity of an OC-like scheme and can be implemented very compactly.

Because of its simplicity, the \texttt{ocUpdate} cannot be expected to be as robust as more general optimization routines, such as the MMA \citep{svanberg_87a} or GCMMA \citep{svanberg_02a}. Therefore, we decided to keep this routine separate from the main code, such that the user can easily switch to other, more general and established optimizers as needed. Details on how to call the MMA within the main code are given in \autoref{App:ksProperties&OCupdate}.

\begin{figure}[t]
 \centering
  \includegraphics[scale = 0.45, keepaspectratio]
   {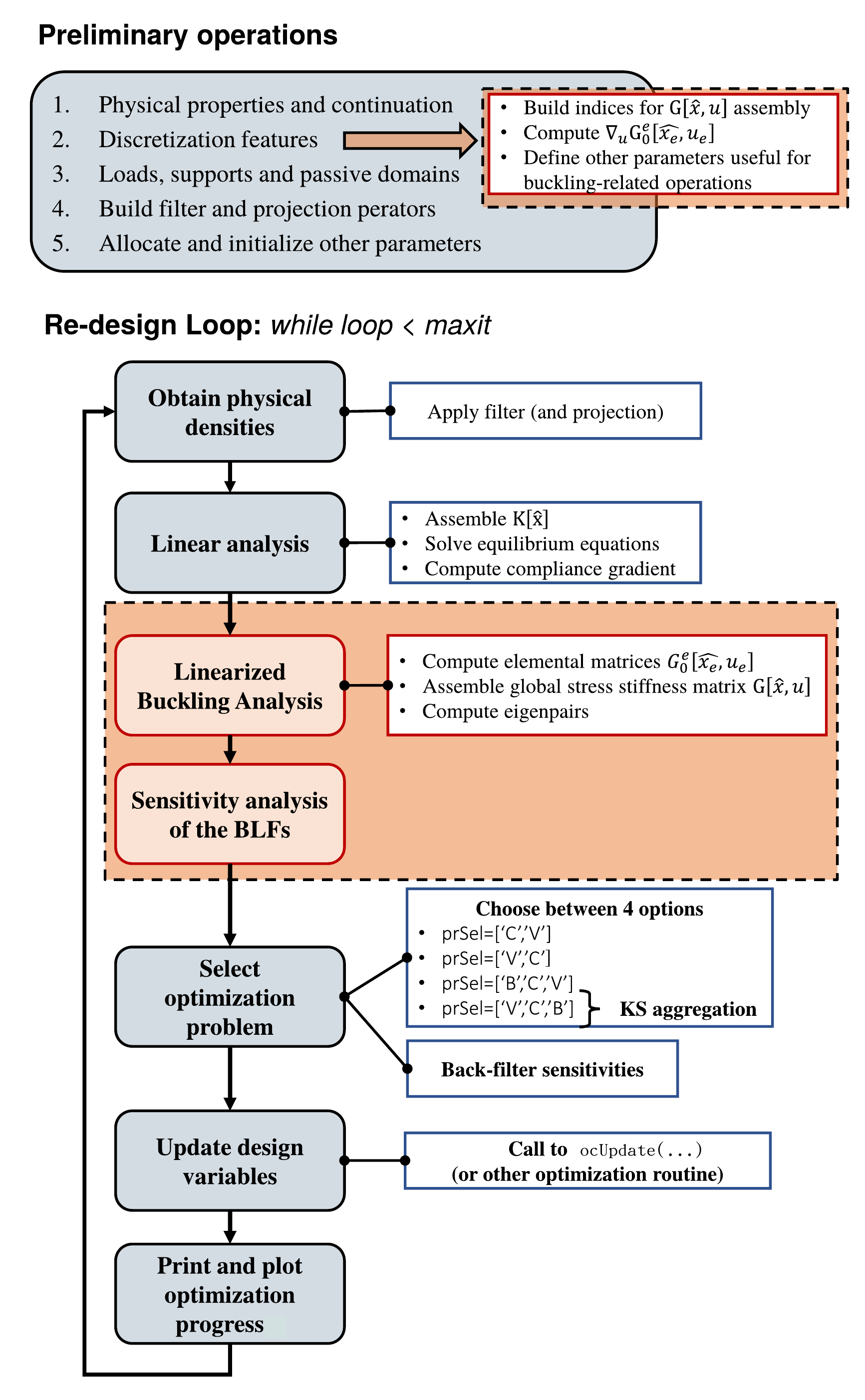}
 \caption{\small{Flowchart for the \texttt{topBuck250} Matlab code. All the operations related to buckling (highlighted with an orange background) are performed only if \texttt{'B'} is among the argument of \texttt{prSel}, otherwise the code essentially reduces to the \texttt{top99neo}. The design variables are updated by the routine \texttt{ocUpdate}}}
 \label{fig:flowchart}
\end{figure}

\section{General structure of the code}
 \label{Sec:codeStructure}

The code is called with the following input data
\begin{lstlisting}[basicstyle=\scriptsize\ttfamily,breaklines=true,numbers=none,frame=single]
topBuck250(nelx,nely,penalK,rmin,ft,ftBC,eta,beta,maxit,ocPar,Lx,penalG,nEig,pAgg,prSel,x0)
\end{lstlisting}
where the first 9 arguments have the same meaning as in \texttt{top99neo} \citep{ferrari-sigmund_20b}. \texttt{ocPar} contains the parameters governing the design update routine, and will be discussed in \autoref{App:ksProperties&OCupdate}.

We introduce the physical length of the domain \texttt{Lx} and the height is obtained by the aspect ratio (line 11). \texttt{penalG} is the penalization factor used for the stress stiffness interpolation \eqref{eq:stiffnessAndStressInterpolations}, \texttt{nEig} is the number of BLFs included in the optimization and \texttt{pAgg} is the initial value of the KS aggregation factor.

\texttt{prSel} is a data structure specifying the optimization problem and the corresponding constraints bounds. This can take the forms
\begin{lstlisting}[basicstyle=\scriptsize\ttfamily,breaklines=true,numbers=none,frame=none]
 prSel={['B','C','V'],[Cmax,Vfmax]};
 prSel={['V','C','B'],[Cmax,Lmin]};
 prSel={['C','V'],Vfmax};
 prSel={['V','C'],Cmax};
\end{lstlisting}
The first substructure \texttt{prSel\{1\}} is a list of 2 or 3 characters, used to select the optimization problem to be solved. The first character defines the objective function and the other(s) the contraint(s), among the criteria introduced in \autoref{Sec:Formulations}. By convention, we set \texttt{'B'} for the BLF, \texttt{'C'} for compliance and \texttt{'V'} for volume fraction. If \texttt{'B'} is not among the set of arguments, all the buckling-related operations are skipped and we essentially recover the \texttt{topo99neo} code, except for the \texttt{ocUpdate} routine (see \autoref{fig:flowchart}). In this case, volume-constrained minimum compliance (\texttt{prSel\{1\}=['C','V']}) and compliance-constrained minimum volume (\texttt{prSel\{1\}=['V','C']}) TO problems can be solved with essentially the same efficiency as in \texttt{top99neo}. Problems \eqref{eq:optProblemMaxBLF} and \eqref{eq:optProblemMinVolcstrBLFandC} are selected by setting \texttt{prSel\{1\} = ['B','C','V']} and \texttt{prSel\{1\} = ['V','B','C']}, respectively.

The second substructure \texttt{prSel\{2\}} contains 1 or 2 numerical values, specifying the constraint(s) bound(s) $\bar{c}$, $\bar{f}$, $\underline{\lambda}$. We assume the following: (1) for volume minimization problems (i.e, if \texttt{prSel\{1\}(1)='V'}) the initial volume fraction is set to $f = 1$ (see line 9); (2) for volume constrained problems, the maximum volume fraction is specified by \texttt{prSel\{2\}(end)=Vfmax}; (3) the compliance upper bound is expressed as a scaling of the initial compliance (e.g., \texttt{Cmax=2.5} means $\bar{c} = 2.5 (\mathbf{F}^{T}\mathbf{u})\mid_{\rm loop = 1}$).

The last argument is a string with the name of a data file (e.g. \texttt{'myData.mat'}), and can be used to specify a non-uniform initial material distribution. This is useful for solving reinforcement problems starting from an initial design (see \autoref{sSec:ssBeam}). We assume that the material distribution is saved as \texttt{"xInitial"} in the data file, and then is assigned to the design variables on line 91 in the code. If this last argument is not passed to \texttt{topBuck250}, the design variables vector is initialized to the uniform material distribution fulfilling the specified volume fraction \texttt{volfrac} (see lines 93-94).

The \texttt{topBuck250} routine starts with some preliminary operations subdivided in the following blocks
\begin{lstlisting}[basicstyle=\scriptsize\ttfamily,breaklines=true,numbers=none,frame=none]
 PRE.1) MATERIAL AND CONTINUATION PARAMETERS
 PRE.2) DISCRETIZATION FEATURES
 PRE.3) LOADS, SUPPORTS AND PASSIVE DOMAINS
 PRE.4) PREPARE FILTER AND PROJECTION OPERATORS
 PRE.5) ALLOCATE AND INITIALIZE OTHER PARAMETERS
\end{lstlisting}

We first define the elastic constants, the continuation schemes for the penalization exponents (\texttt{penalK}, \texttt{penalG}), for the projection parameter \texttt{beta} and for the aggregation exponent \texttt{pAgg} (lines 3-7). We remark that, as in \texttt{top99neo}, the continuation scheme on a given parameter is specified as
\begin{lstlisting}[basicstyle=\scriptsize\ttfamily,breaklines=true,numbers=none,frame=none]
 parCont={istart,maxPar,isteps,deltaPar};
\end{lstlisting}
such that when \texttt{loop>=istart}, the parameter is increased by \texttt{deltaPar} every \texttt{isteps}, up to the value \texttt{maxPar}. The continuation is applied at lines 245-246, by the function defined on line 8. The check on lines 247-248 keeps track of the steps where continuation is applied, and allows the user to restart the asymptotes history (see \autoref{App:ksProperties&OCupdate}).

The discretization features are built between lines 11-59. Up to line 32 the operations are the same as in \texttt{top99neo}, and we remark that the indices \texttt{iK} and \texttt{jK}, used for the $K$ matrix assembly (lines 27-32), refer to the lower-half of the symmetric matrix and are built as explained in \cite{ferrari-sigmund_20b}.

Lines 33-62 are executed only if buckling is among the optimization criteria. We introduce the non dimensional elasticity matrix \texttt{Cmat0}, the physical dimensions of the elements \texttt{xe}, and the Gauss nodes (\texttt{xiG,etaG}) and weights (\texttt{wxi,weta}) adopted for numerical quadratures (lines 34-36). Some operators useful in order to compactly perform the stress analysis are introduced in lines 37-39 (see \autoref{App-notesDiscretizationOperators} for details). The indices \texttt{iG} and \texttt{jG} used for the stress stiffness matrix assembly, the derivative of this latter with respect to the displacement vector, and other auxiliary variables are defined between lines 40-59. These operations are all precisely discussed in \autoref{Sec:BucklingFast}.

The KS function and its derivative, as specified by \eqref{eq:ksAggregationLambda1} and \eqref{eq:ksAggregationSens} are defined as functions on lines 60 and 61. The specification of loads, boundary conditions and passive elements follows the same concept of \texttt{top88} and \texttt{top99neo}, and the instructions given by default on lines 64-69 correspond to the compressed column example of \autoref{sSec:ssBeam}.

The filter operator is built in \texttt{PRE.4)}, using Dirichlect (\texttt{bcF='D'}) or Neumann (\texttt{bcF='N'}) boundary conditions, and lines 78-83 define the $\eta$ projection (\texttt{prj}) and its $\eta$- and $\tilde{x}_{e}$-derivatives (\texttt{deta}, \texttt{dprj}). Block \texttt{PRE.5)} is basically as in \texttt{top99neo}. The only notable difference is that lines 90-95 allow the user to initialize the design variables to the initial guess specified by \texttt{x0}, or to the uniform material distribution.

The redesign loop starts at line 98 and consists of the following blocks
\begin{lstlisting}[basicstyle=\scriptsize\ttfamily,breaklines=true,numbers=none,frame=none]
 RL.1) COMPUTE PHYSICAL DENSITY FIELD
 RL.2) SETUP AND SOLVE EQUILIBRIUM EQUATIONS
 RL.3) BUILD STRESS STIFFNESS MATRIX
 RL.4) SOLVE BUCKLING EIGENVALUE PROBLEM
 RL.5) SENSITIVITY ANALYSIS OF BLFs
 RL.6) SELECT OBJECTIVE FUNCTION AND CONSTRAINTS
 RL.7) UPDATE DESIGN VARIABLES
 RL.8) PRINT AND PLOT RESULTS
\end{lstlisting}

\texttt{RL.1)} and \texttt{RL.2)} perform the operations needed for computing the linear compliance and its derivative. The operations in blocks \texttt{RL.3)} to \texttt{RL.5)} are discussed in detail in \autoref{Sec:BucklingFast} and are executed only if buckling is among the optimization criteria (see \autoref{fig:flowchart}). Such a modular organization of the code allows the user to skip the most computationally intensive parts if the purpose is to solve standard compliance or volume minimization problems, e.g in order to generate an initial design to be then reinforced against buckling.

Up to this point the script has the greatest generality and, depending on the problem selected, provides the function values and gradients for the three physical responses: compliance, volume fraction and $\mu_{i}$ parameters, that are related to BLFs. Then, \texttt{RL.6)} and \texttt{RL.7)} are customized for the particular optimization problems introduced in \autoref{sSec:SolutionScheme} and for performing the design update with the \texttt{ocUpdate} routine, as described in \autoref{App:ksProperties&OCupdate}. However, we remark that the extension of the code for considering other combinations of the three response criteria, as well as for using more general optimization routines, requires the user to add/modify only a few lines in \texttt{RL.6)} and \texttt{RL.7)}.

The last block (\texttt{RL.8)}) contains some commands for print and plot. If buckling is selected, lines 231-238 generate three sub-plots. The top one shows the current topology, the bottom left one shows the evolution of the constraint functions, together with their KS aggregation, and the evolution of the objective (volume fraction or $J^{KS}[\mu_{i}]^{-1}$). The bottom right plot shows the evolution of the four lowest BLFs, respectively. Otherwise, if buckling is not selected, a single plot showing the evolution of the design is generated.

\section{Setup of the buckling eigenvalue problem and sensitivity analysis of buckling load factors}
 \label{Sec:BucklingFast}

In order to setup the eigenvalue buckling analysis \eqref{eq:eigenvalueEquation} we need to assemble the stress stiffness matrix $G[\hat{\mathbf{x}},\mathbf{u}]$ from the elemental ones $G^{e}[\hat{x}_{e},\mathbf{u}_{e}]$, $e = 1\ldots, m$.

The elemental matrices cannot be expressed as a scaling of a constant matrix (like $K^{e}=E_{K}(\hat{x}_{e})K^{e}_{0}$), but must be evaluated for each element, as they depend on the local displacement field. Within a scripting language such as Matlab, where loops and memory allocation considerably affect the performance, this becomes very expensive even for medium-size problems. However, a careful inspection of the structure of $G^{e}[\hat{x}_{e},\mathbf{u}_{e}]$ allows a very efficient implementation, requiring cheap vectorized matrix products only.

The first step is to compute the stress components $\sigma^{e}_{x}$, $\sigma^{e}_{y}$, $\tau^{e}_{xy}$ at the element level, and arrange them into the matrix $\boldsymbol{\sigma}^{e}_{0}$. For a $k$-noded isoparametric element \citep{book:zienkiewicz-taylor06} this is done through the following relationship
\begin{equation}
 \label{eq:computeStresses}
  \boldsymbol{\sigma}^{e}_{0} = D B_{0} \mathbf{u}^{e}
\end{equation}
where $D$ is the elasticity matrix and $B_{0}$ is the linearized strain-displacement pseudo-differential operator. 

Using the operators defined in \autoref{App-notesDiscretizationOperators}, \autoref{eq:computeStresses} is evaluated  simultaneously for all the elements by the following two lines
\begin{lstlisting}[basicstyle=\scriptsize\ttfamily,breaklines=true,numbers=none,frame=single]
 % ------ stress at current Gauss point (nElemx3)
 gradN = (dN(0,0)*xe)\dN(0,0);
 sGP = (Cmat0*Bmat(gradN)*U(cMat)')';
\end{lstlisting}
where the element centroid is the stress super-convergent point for the $\mathcal{Q}_{4}$ bilinear element \citep{book:wahlbin1995}.

\subsection{Setup of the stress stiffness matrix}
 \label{sSec:fastKsigma}

The stress stiffness matrix stems from the following contribution in the linearization of the virtual work equation \citep{book:zienkiewicz-taylor06}
\begin{equation}
 \label{eq:weakformLB}
  G^{e}[\mathbf{u}_{e}] :=
  \int\limits_{\Omega_{0}} B^{T}_{1}
  \mathsf{T}[\boldsymbol{\sigma}^{e}_{0}] B_{1} \: {\rm d}\Omega_{0}
\end{equation}
where $B_{1}$ discretizes the deformation gradient and for a $\mathcal{Q}_{4}$ bilinear element reads \citep{book:crisfield91}
\begin{equation*}
 \label{eq:discretizedFgrad}
  B_{1} = 
  \left[
  \begin{array}{cccccccc}
    \partial_{x}N_{1} & 0 & \partial_{x}N_{2} & 0 & \partial_{x}N_{3} & 0 & \partial_{x}N_{4} & 0 \\
    \partial_{y}N_{1} & 0 & \partial_{y}N_{2} & 0 & \partial_{y}N_{3} & 0 & \partial_{y}N_{4} & 0 \\
    0 & \partial_{x}N_{1} & 0 & \partial_{x}N_{2} & 0 & \partial_{x}N_{3} & 0 & \partial_{x}N_{4} \\
    0 & \partial_{y}N_{1} & 0 & \partial_{y}N_{2} & 0 & \partial_{y}N_{3} & 0 & \partial_{y}N_{4}
  \end{array}
  \right]
\end{equation*}
and the stress components are arranged as
\begin{equation}
 \mathsf{T}[\boldsymbol{\sigma}^{e}_{0}] = 
 I_{2} \otimes \boldsymbol{\sigma}^{e}_{0} = 
 \left[
  \begin{array}{cccccccc}
   \sigma_{x} & \tau_{xy} & 0 & 0 \\
   \tau_{xy} & \sigma_{y} & 0 & 0 \\
   0 & 0 & \sigma_{x} & \tau_{xy} \\
   0 & 0 & \tau_{xy} & \sigma_{y}
  \end{array}
 \right]
\end{equation}
where $I_{2}$ is the identity matrix of order 2.

Expanding the product $B^{T}_{1}\mathsf{T}[\boldsymbol{\sigma}^{e}_{0}] B_{1}$ we obtain the following structure for the integrand of \eqref{eq:weakformLB} \citep{book:deborst2012}
\begin{equation}
 \label{eq:elementGeometricMatrixFull}
 G^{e}_{\left\{\mathcal{Q}_{4}\right\}} = 
 \left[
  \begin{array}{ccccccccc}
   \textcolor{red}{z_{11}} & 0 & z_{12} & 0 & z_{13} & 0 & z_{14} & 0 \\
   0 & \textcolor{cyan}{z_{11}} & 0 & z_{12} & 0 & z_{13} & 0 & z_{14} \\
   \textcolor{red}{z_{21}} & 0 & \textcolor{red}{z_{22}} & 0 & z_{23} & 0 & z_{24} & 0 \\
   0 & \textcolor{cyan}{z_{21}} & 0 & \textcolor{cyan}{z_{22}} & 0 & z_{23} & 0 & z_{24} \\
   \textcolor{red}{z_{31}} & 0 & \textcolor{red}{z_{32}} & 0 & \textcolor{red}{z_{33}} & 0 & z_{34} & 0 \\
   0 & \textcolor{cyan}{z_{31}} & 0 & \textcolor{cyan}{z_{32}} & 0 & \textcolor{cyan}{z_{33}} & 0 & z_{34} \\
   \textcolor{red}{z_{41}} & 0 & \textcolor{red}{z_{42}} & 0 & \textcolor{red}{z_{43}} & 0 & \textcolor{red}{z_{44}} & 0 \\
   0 & \textcolor{cyan}{z_{41}} & 0 & \textcolor{cyan}{z_{42}} & 0 & \textcolor{cyan}{z_{43}} & 0 & \textcolor{cyan}{z_{44}}
  \end{array}
 \right]
\end{equation}
and we identify only 10 independent coefficients out of the 64,  highlighted in red in \eqref{eq:elementGeometricMatrixFull}. Moreover, the generic $z_{ik}$ ($i, k = 1, \ldots, 4$) can be expressed as a linear combination of the stress components
\begin{equation}
 \label{eq:genericStressMatrixCoefficient}
  \begin{aligned}
   z_{ik} & = \sigma_{x} a_{k}a_{i} + \sigma_{y} b_{k}b_{i} + \tau_{xy} \left( b_{k}a_{i} + a_{k}b_{i} \right)
  \end{aligned}
\end{equation}
where we have set $a_{i} = \partial_{x}N_{i}$ and $b_{i} = \partial_{y}N_{i}$.

\begin{table}[t]
 \centering
  \begin{tabular}{p{0.75cm}|>{\centering}p{0.3cm}>{\centering}p{0.3cm}>{\centering}p{0.3cm}>{\centering}p{0.3cm}>{\centering}p{0.3cm}>{\centering}p{0.3cm}>{\centering}p{0.3cm}>{\centering}p{0.3cm}>{\centering}p{0.3cm}p{0.3cm}}
    & 1 & 2 & 3 & 4 & 5 & 6 & 7 & 8 & 9 & 10 \\
    $\mathcal{M}_{\rm odd}$  & (1,1) & (3,1) & (5,1) & (7,1) & (3,3) & (5,3) & (7,3) & (5,5) & (7,5) & (7,7) \\ 
    $\mathcal{M}_{\rm even}$ & (2,2) & (4,2) & (6,2) & (8,2) & (4,4) & (6,4) & (8,4) & (6,6) & (8,6) & (8,8)
  \end{tabular}
 \caption{\small{Mapping $\mathcal{M}_{\rm odd}:\{ 1 : 10 \} \mapsto (i,j)$ for the 10 unique coefficients of $\textsf{Z}$ and the \emph{lower} symmetric part of $G^{e}$. $\mathcal{M}_{\rm even}$ in the second row is for even columns}}
 \label{tab:endexingTenCoeffge}
\end{table}

A compact representation of the 10 independent coefficients in \eqref{eq:elementGeometricMatrixFull} is given by the array $\mathsf{Z}\in\mathbb{R}^{m\times 10}$, which is obtained through the Hadamard matrix product $\textsf{Z} = \textsf{B}^{T}\odot \textsf{S}$. Here, $\mathbb{R}^{m\times 3} \ni \textsf{S} = \left[ \sigma^{e}_{x}, \sigma^{e}_{y}, \tau^{e}_{xy} \right]_{e=1:m}$ collects the stress components for all the elements, and we have introduced the array $\mathsf{B}_{[3 \times 10]}$
\begin{equation}
 \textsf{B} = \left[
  \begin{array}{c}
   a_{\ell(i, 1)}a_{\ell(i, 2)} \\
   b_{\ell(i, 1)}b_{\ell(i, 2)} \\
   b_{\ell(i, 2)}a_{\ell(i, 1)} +
   a_{\ell(i, 2)}b_{\ell(i, 1)}
  \end{array}
 \right]^{i=1:10}
\end{equation}
where $\ell\in\mathbb{N}^{10\times 2}$ is the set of indices mapping each $\mathsf{Z}_{ei}$ ($i=1, \ldots 10$) to the corresponding $z_{\ell(i,1)\ell(i,2)}$ within the pattern of non-zero elements of \eqref{eq:elementGeometricMatrixFull}, i.e.
\[
  \ell =
  \left[
  \begin{array}{cccccccccc}
   1 & 2 & 3 & 4 & 2 & 3 & 4 & 3 & 4 & 4 \\
   1 & 1 & 1 & 1 & 2 & 2 & 2 & 3 & 3 & 4
  \end{array}
 \right]^{T}
\]

The columns of $\textsf{Z}$ collect the coefficients highlighted in the reduced pattern here below, sorted column-wise
\begin{equation}
 \label{eq:elementGeometricMatrixReducedPattern}
 G^{e,\rm unique}_{\left\{\mathcal{Q}_{4}\right\}} = 
 \left[
  \begin{array}{cccc}
   \textcolor{blue}{z_{11}} & z_{12} & z_{13} & z_{14} \\
   \textcolor{blue}{z_{21}} & \textcolor{blue}{z_{22}} & z_{23} & z_{24} \\
   \textcolor{blue}{z_{31}} & \textcolor{blue}{z_{32}} & \textcolor{blue}{z_{33}} & z_{34} \\
   \textcolor{blue}{z_{41}} & \textcolor{blue}{z_{42}} & \textcolor{blue}{z_{43}} & \textcolor{blue}{z_{44}} \\
  \end{array}
 \right]
\end{equation}

Thus, we obtain the independent stress stiffness matrix coefficients for all the elements at once, and we just need a small loop to evaluate $\textsf{Z}$ on the Gauss points. We remark that $\textsf{Z}$ actually contains non-dimensional strain combinations, and $K^{e}_{0}$ and $G^{e}_{0}$ are scaled by the interpolations $E_{K}\left( \hat{x}_{e} \right)$ and $E_{G}\left( \hat{x}_{e} \right)$, to obtain the dimensional matrices $K^{e}$ and $G^{e}$, right before the global assembly. The whole procedure for arriving at $G[\hat{\mathbf{x}},\mathbf{u}]$ is implemented between lines 124-148.

\begin{figure*}
 \centering
  \subfloat[2D discretization]{
   \includegraphics[scale = 0.36, keepaspectratio]
    {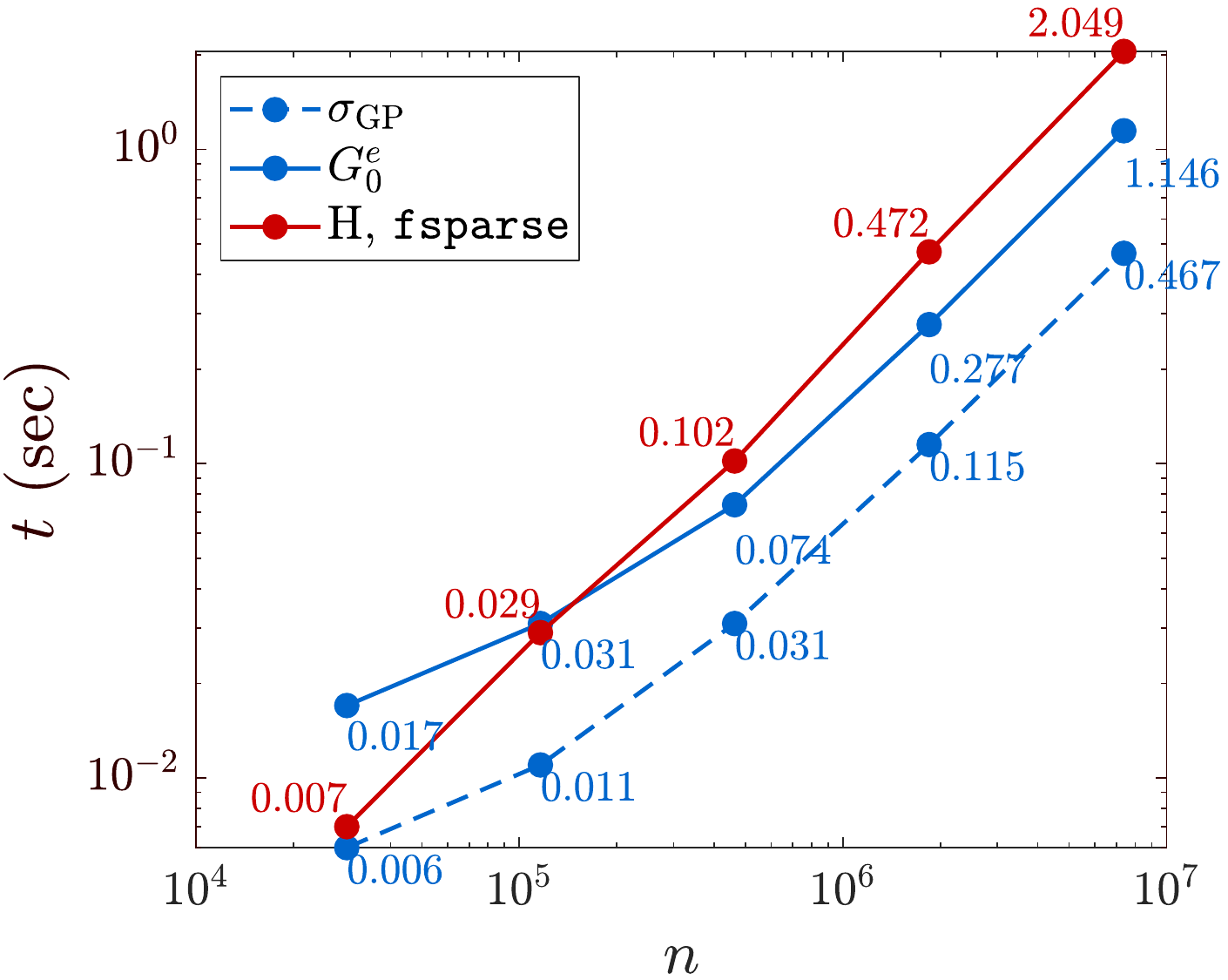}} \quad
  \subfloat[3D discretization]{
   \includegraphics[scale = 0.36, keepaspectratio]
    {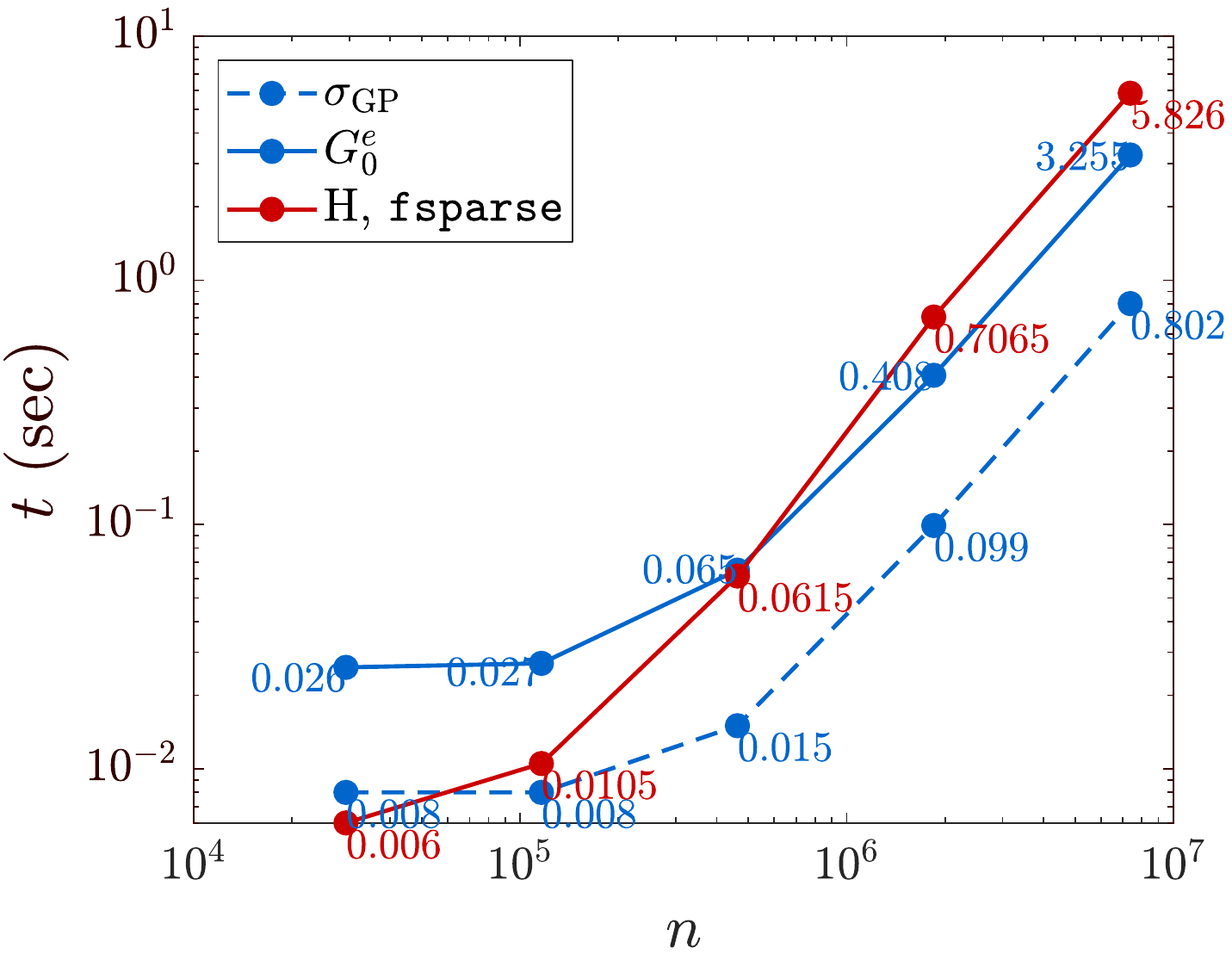}} \quad
  \subfloat[]{
   \includegraphics[scale = 0.36, keepaspectratio]
    {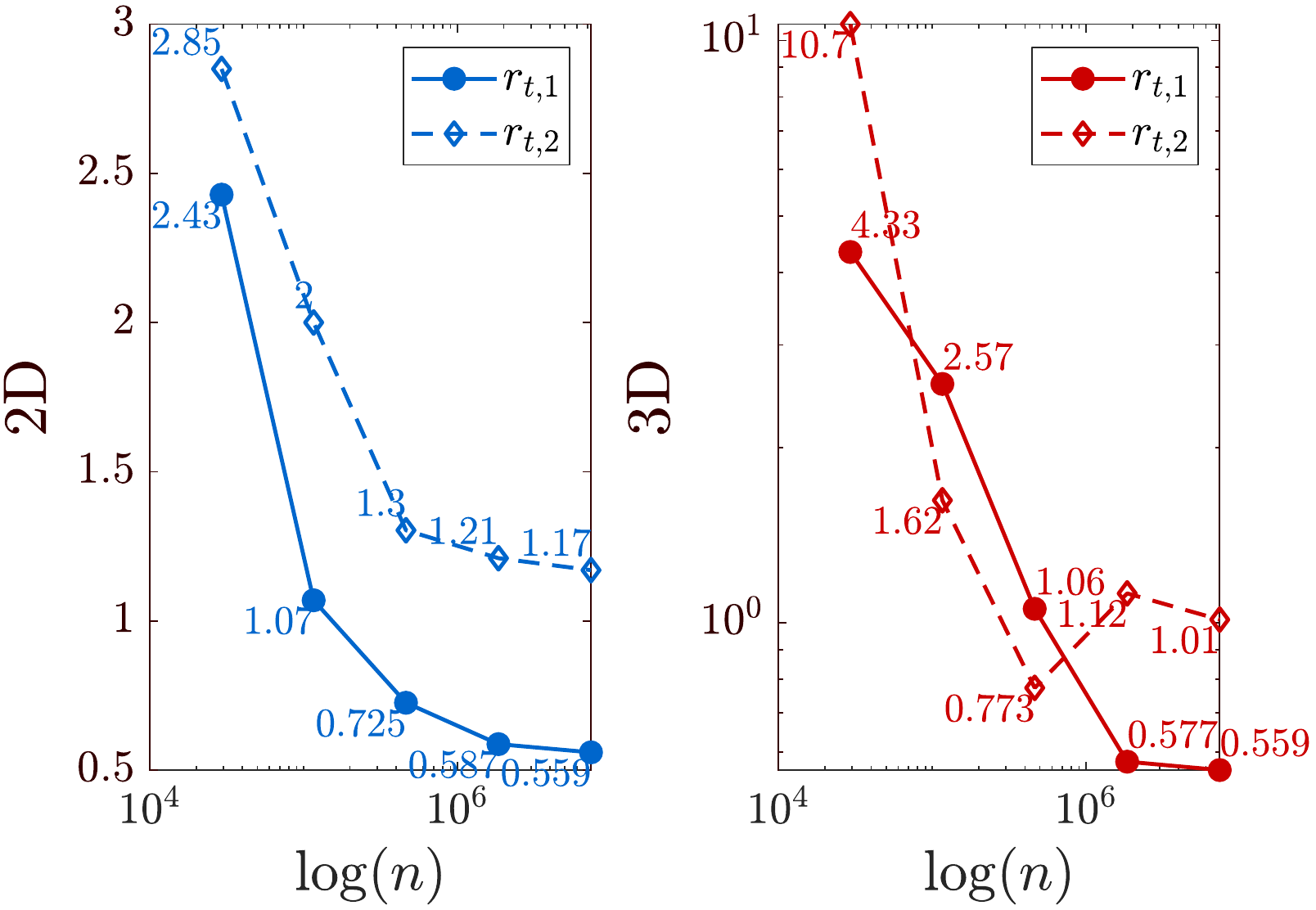}}
 \caption{\small{Scaling of CPU time for the construction of the stress stiffness matrix $G$ (a,b). The overall time is split between that for the stresses computation ($\sigma_{GP}$}), for computing the compact representation $\textsf{Z}$ of the elemental matrices $(G^{e}_{0})$, and that for assembling the global matrix ({\rm H},\texttt{fsparse}). In (c) we see the scaling of the ratios $r_{t, 1} = \frac{\sigma_{GP}+G^{e}_{0}}{{\rm H},\texttt{fsparse}}$, measuring the efficiency of the elemental matrices setup, and $r_{t, 2}$, comparing the overall time for setting up $G$ with that for setting up $K$}
 \label{fig:ksSetup}
\end{figure*}

The global stress stiffness matrix $G$ can be assembled as discussed by \cite{ferrari-sigmund_20b}, upon considering the mapping $\mathcal{M}_{\rm odd}$ between the 10 columns of $\textsf{Z}$ and the \emph{odd} columns of the \emph{lower} symmetric part of $G^{e}$ (see \autoref{tab:endexingTenCoeffge}). Mapping to \emph{even} columns is then easily obtained by noting that $(i,j)_{\rm even} = (i,j)_{\rm odd} + (1,1)$, and gives the coefficients highlighted in cyan in \eqref{eq:elementGeometricMatrixFull}. The indices \texttt{iG} and \texttt{jG} are subsets of \texttt{iK} and \texttt{jK}. Considering the column-wise sorting of the \emph{lower symmetric part} of  \eqref{eq:elementGeometricMatrixFull}, it is immediate to recognize that $\mathcal{M}_{\rm odd}$ targets the elements $\{1,3,5,7,16,18,20,27,29,34\}$.

Thus, the assembly indices for both $K$ and $G$ are built by the following instructions, (see lines 26-42)
\begin{lstlisting}[basicstyle=\scriptsize\ttfamily,breaklines=true,numbers=none,frame=single]
 [sI,sII]=deal([]);                         % Line 26
 for j=1:8                                  % Line 27
    sI = cat(2,sI,j:8);                     % Line 28
    sII = cat(2,sII,repmat(j,1,8-j+1));     % Line 29
 end                                        % Line 30
 [iK,jK] = deal(cMat(:,sI)', cMat(:,sII)'); % Line 31
 Iar = sort([iK(:),jK(:)],2,'descend');     % Line 32
 indM = [1,3,5,7,16,18,20,27,29,34];        % Line 40
 [iG,jG] = deal(iK(indM,:),jK(indM,:));     % Line 41
 IkG = sort([iG(:),jG(:)],2,'descend');     % Line 42
\end{lstlisting}

Finally, the lower symmetric part of $G$ is assembled by using \texttt{fsparse} \citep{engblom-lukarski_16a} (see lines 146-147)
\begin{lstlisting}[basicstyle=\scriptsize\ttfamily,breaklines=true,numbers=none,frame=single]
G=fsparse(IkG(:,1)+1,IkG(:,2)+1,sG,[nDof,nDof])+
     fsparse(IkG(:,1),IkG(:,2),sG,[nDof,nDof]);
\end{lstlisting}

The same procedure can be extended to a 3D discretization with $\mathcal{H}_{8}$ trilinear elements. \autoref{eq:genericStressMatrixCoefficient} becomes ($\left( i, k \right) = 1, \ldots, 8$)
\begin{equation*}
 \label{eq:GenericStressMatrixCoefficientH8}
  \begin{aligned}
   z_{ik} & = \sigma_{x} a_{k}a_{i} + \sigma_{y} b_{k}b_{i} + \sigma_{z} c_{k}c_{i} + \tau_{xy} \left( b_{k}a_{i} + a_{k}b_{i} \right) \\
   & + \tau_{xz} \left( a_{k}c_{i} + c_{k}a_{i} \right) + \tau_{yz} \left( b_{k}c_{i} + c_{k}b_{i} \right)
  \end{aligned}
\end{equation*}
where $c_{i} = \partial_{z}N_{i}$ and the pattern of $G^{e}_{\left\{\mathcal{H}_{8}\right\}}$ is formally similar to that of \eqref{eq:elementGeometricMatrixFull}. For the 3D case, computational cuts are much more substantial, as there are just 36 independent coefficients out of 576.

We refer to \autoref{fig:ksSetup} for considerations about the efficiency. Plots (a,b) show the scaling curves for the CPU time spent on computing the elemental stresses ($t[\sigma_{GP}]$) and the elemental geometric stiffness matrices $t[G^{e}_{0}]$. As $n$ increases, the sum of these two times becomes smaller than the time spent on the $G$ assembly operation ($t[{\rm H},\texttt{fsparse}]$). We notice that the assembly is a bit cheaper for the stress stiffness matrix than for the stiffness one, due to the larger sparsity of $G$ (only 1/2 of the $z_{ik} \neq 0$, whereas $K^{e}_{0}$ is full).

\autoref{fig:ksSetup} (c) shows the trend of two coefficients measuring the efficiency of the operations for building the global stress stiffness matrix $G$. $r_{t, 1}$ is the ratio between $t[\sigma_{GP}] + t[G^{e}_{0}]$ (i.e. the overall time for computing the independent coefficients of the elemental stress stiffness matrices) and $t[{\rm H},\texttt{fsparse}]$. For small discretizations we have $r_{t, 1}>1$. This is reasonable because the assembly operation is extremely cheap and some operations for computing $G^{e}_{0}$ (e.g. the short loop on Gauss knots) are not amortized by the efficiency of the vector products. However, as $n$ increases $r_{t,1}$ drops as the computation of $G^{e}_{0}$ takes advantage of the fully vectorized operations. $r_{t, 2}$ is the ratio between the overall times for setting up $G$ and that for $K$. Also this term decreases as the DOFs number $n$ becomes larger, and for very large 2D discretizations the computation of $G$ is about $17\%$ in 2D, and only $1\%$ in 3D, more expensive than computing $K$.

Therefore, the procedure we have described cuts the cost for setting up the buckling eigenvalue problem to the same as the one for setting up the linear analysis.

\subsection{Sensitivity expressions}
 \label{sSec:sensitivityBuckling}

The $\hat{x}_{e}$-derivatives of the elemental matrices, needed for the expressions \eqref{eq:sensitivityComplianceVolume} and \eqref{eq:sensitivityMu}, are given by
\begin{equation}
 \label{eq:DerivativeOfLocalMatrices}
 \begin{aligned}
  \frac{\partial K^{e}}{\partial \hat{x}_{e}} = \frac{\partial E_{K}\left( \hat{x}_{e} \right)}
  {{\partial} \hat{x}_{e}} K^{e}_{0}
  \ , \qquad
  \frac{\partial G^{e}}{\partial \hat{x}_{e}} = \frac{\partial E_{G}\left( \hat{x}_{e} \right)}
  {{\partial} \hat{x}_{e}} G^{e}_{0}
 \end{aligned}
\end{equation}
and the gradient of the compliance is computed as usual \citep{andreassen-etal_11a,ferrari-sigmund_20b}, expressing the local product $\mathbf{u}^{T}_{e}K^{e}_{0}\mathbf{u}_{e}$ through the connectivity matrix (line 121).

Below we give details about how to take advantage of the compact representation introduced for $G^{e}_{0}$ (i.e. $\mathsf{Z}\in\mathbb{R}^{m\times 10}$) when expressing each term of \eqref{eq:sensitivityMu}.

\paragraph{Contribution 1.} To express this at the element level we need to account for the structure of $\textsf{Z}$. Writing out explicitly the double product we obtain
\begin{small}
\begin{equation}
 \label{eq:doubleProductvGv}
 \begin{aligned}
  \boldsymbol{\varphi}^{T}G^{e}_{0}\boldsymbol{\varphi} & = 
      \left( \varphi^{2}_{1} + \varphi^{2}_{2} \right) z_{11} + \left( \varphi^{2}_{3} + \varphi^{2}_{4} \right) z_{22} \\
  & + \left( \varphi^{2}_{5} + \varphi^{2}_{6} \right) z_{33} + \left( \varphi^{2}_{7} + \varphi^{2}_{8} \right) z_{44} \\
  & + 2\left( \varphi_{1}\varphi_{3} + \varphi_{2}\varphi_{4} \right) z_{21} + 2\left( \varphi_{1}\varphi_{5} + \varphi_{2}\varphi_{6} \right) z_{31} \\
  & + 2(\varphi_{1}\varphi_{7} + \varphi_{2}\varphi_{8}) z_{41} + 2(\varphi_{3}\varphi_{5} + \varphi_{4}\varphi_{6}) z_{32} \\
  & + 2(\varphi_{3}\varphi_{7} + \varphi_{4}\varphi_{8}) z_{42} + 2(\varphi_{5}\varphi_{7} + \varphi_{6}\varphi_{8}) z_{43}
 \end{aligned}
\end{equation}
\end{small}
where we have exploited the symmetry of $z_{ij}$ and we have set $\boldsymbol{\varphi} = \boldsymbol{\varphi}_{e}$ for simplicity. The same result of \eqref{eq:doubleProductvGv} can be obtained with the Hadamard product $\tilde{\textsf{Z}}\odot\textsf{p}$, where $\mathsf{p}\in\mathbb{R}^{10\times m}$ collects the components of $\boldsymbol{\varphi}$ as shown in \autoref{tab:IndexingPforModes}, for all the elements. A quick glance at the subscripts in \autoref{tab:IndexingPforModes} reveals that the array $\mathsf{p}$ can be built by using the indices \texttt{iG} and \texttt{jG}. Indeed, we recognize $\mathsf{p}_{i} = \varphi_{\mathcal{M}_{\rm odd}(i, 2)}\varphi_{\mathcal{M}_{\rm odd}(i, 1)} + \varphi_{\mathcal{M}_{\rm even}(i, 2)}\varphi_{\mathcal{M}_{\rm even}(i, 1)}$; thus we formally express $\mathsf{p} = \boldsymbol{\varphi}(\tilde{\mathcal{M}})$, where $\tilde{\mathcal{M}}\in\mathbb{N}^{10\times m}$.

To have consistency with \eqref{eq:doubleProductvGv} the elements on columns $\{2,3,4,6,7,9\}$ of the array $\textsf{Z}$ must be doubled, obtaining $\tilde{\textsf{Z}}$. This first contribution could be computed by the following operations
\begin{lstlisting}[basicstyle=\scriptsize\ttfamily,breaklines=true,numbers=none,frame=single]
t2ind=[2,3,4,6,7,9];            % defined on Line 40
a1=reshape(IkG(:,2),10,nEl)';   % defined on Line 43
a2=reshape(IkG(:,1),10,nEl)';   % defined on Line 43
dkeG=dsG.*Z;                    % x-derivative of Ge
dkeG(:,t2ind)=2*dkeG(:,t2ind);  % 2x t2ind columns
for j=1:nEig
   t=phi(:,j);
   vC=t(a1).*t(a2)+t(a1+1).*t(a2+1);
   phiDGphi(:,j)=sum(dkeG.*vC,2);
end
\end{lstlisting}
similar to those in the code of \autoref{App:matlabCodeBuckling}.

\begin{table}[t]
 \centering
  \begin{tabular}{c|c}
    $\textsf{p}_{1}$ & $(\varphi^{2}_{1} + \varphi^{2}_{2})$ \\
    $\textsf{p}_{2}$ & $(\varphi_{1}\varphi_{3} + \varphi_{2}\varphi_{4})$ \\
    
    $\textsf{p}_{3}$ & $(\varphi_{1}\varphi_{5} + \varphi_{2}\varphi_{6})$ \\
    $\textsf{p}_{4}$ & $(\varphi_{1}\varphi_{7} + \varphi_{2}\varphi_{8})$ \\
    $\textsf{p}_{5}$ & $(\varphi^{2}_{3} + \varphi^{2}_{4})$ \\
    $\textsf{p}_{6}$ & $(\varphi_{3}\varphi_{5} + \varphi_{4}\varphi_{6})$ \\
    $\textsf{p}_{7}$ & $(\varphi_{3}\varphi_{7} + \varphi_{4}\varphi_{8})$ \\
    $\textsf{p}_{8}$ & $(\varphi^{2}_{5} + \varphi^{2}_{6})$ \\
    $\textsf{p}_{9}$ & $(\varphi_{5}\varphi_{7} + \varphi_{6}\varphi_{8})$ \\
    $\textsf{p}_{10}$ & $(\varphi^{2}_{7} + \varphi^{2}_{8})$
  \end{tabular}
 \caption{\small{Mapping of the eigenvector components into the array $\mathsf{p}$ used for writing the product $\boldsymbol{\varphi}^{T}G^{e}_{0}\boldsymbol{\varphi} \equiv \tilde{\mathsf{Z}}\odot\mathsf{p}$}}
 \label{tab:IndexingPforModes}
\end{table}

\paragraph{Contribution 2}

The double product $\boldsymbol{\varphi}^{T}(\partial_{e}K)\boldsymbol{\varphi}$ can be computed at the element level as done for $\partial_{e}c$, upon replacing $\partial_{e}E_{K}(x_{e})$ by $\partial_{e}E_{G}(x_{e})$ and $\mathbf{u}$ by $\boldsymbol{\varphi}$.

\paragraph{Contribution 3}

The derivative of the non-dimensional elemental stress stiffness matrix $G^{e}_{0}$ with respect to each component $u_{i}$ of the displacement vector $\mathbf{u}$ reads (see \autoref{eq:weakformLB})
\begin{equation}
 \label{eq:derivativeGu}
 \begin{aligned}
  \frac{\partial G^{e}_{0}}{\partial u_{i}}(\mathbf{u}) & =
  \frac{\partial}{\partial u_{i}}
  \int\limits_{\Omega_{e}} B^{T}_{1} \left( I_{2} \otimes
  \frac{\partial \boldsymbol{\sigma}^{e}}{\partial u_{i}} \right) B_{1}
  \: {\rm d}\Omega_{e} \\ & =
   \int\limits_{\Omega_{e}} B^{T}_{1} \left( I_{2} \otimes C B_{0}
  \frac{\partial \mathbf{u}_{e}}{\partial u_{i}} \right) B_{1}
  \: {\rm d}\Omega_{e} \\ & =
   \int\limits_{\Omega_{e}} B^{T}_{1} \left( I_{2} \otimes C B_{0}
  \delta_{ie} \right) B_{1} \: {\rm d}\Omega_{e}
 \end{aligned}
\end{equation}
where $\delta_{ie} = 1$ if $u_{i}$ belongs to the displacements $\mathbf{u}_{e}$ on element $e$, and $\delta_{ie} = 0$ otherwise.

Since \eqref{eq:derivativeGu} is independent of $\hat{\mathbf{x}}$, a simple way to obtain the gradient $\nabla_{\mathbf{u}}G^{e}_{0}$ is to pre-compute the operator $\nabla_{\mathbf{u}}\textsf{Z} \in \mathbb{R}^{10\times 8}$. Following the same reduced storage format of $\mathsf{Z}$, the $i$-th column of $\nabla_{\mathbf{u}}\textsf{Z}$ contains the 10 independent coefficients of the stress stiffness matrix corresponding to the displacement vector $\mathbf{u} \in \mathbb{R}^{8\times 1}$, $u_{i} = 1$, $u_{j\neq i} = 0$; namely,
\begin{equation}
 \label{eq:gradientGwrtU}
  \nabla_{\mathbf{u}}\textsf{Z} =
  \left[
  G^{e}_{0}(u_{1}=1),
  G^{e}_{0}(u_{2}=1),
  \ldots,
  G^{e}_{0}(u_{8}=1)
  \right]
\end{equation}

\autoref{eq:gradientGwrtU}, called \texttt{dZdu} in the code, is computed between lines 44-58, using the compact definitions of $B_{0}$ and $B_{1}$ given in \autoref{App-notesDiscretizationOperators}. Then, each adjoint load can be computed by expressing the following double product
\begin{equation}
 \label{eq:adjLoadsImplementation}
 \boldsymbol{\varphi}^{T}_{i}
 (\nabla_{\mathbf{u}}G)\boldsymbol{\varphi}_{i} = 
 (\nabla_{\mathbf{u}}\tilde{\textsf{Z}})\odot\textsf{p}
\end{equation}
as described for Contribution $\circled{1}$, where $\mathsf{p}$ is as in \autoref{tab:IndexingPforModes} and $\nabla_{\mathbf{u}}\tilde{\mathsf{Z}}$ is obtained from $\nabla_{\mathbf{u}}\mathsf{Z}$ by doubling the elements of columns $\{2,3,4,6,7,9\}$. Once the adjoint problem has been solved for $\mathbf{w}_{i}$, the computation of the last sensitivity term only requires another application of step $\circled{2}$, with $\mathbf{u}$ and $\mathbf{w}_{i}$ as left and right vectors, respectively.

Sensitivity calculations are implemented in the code of \autoref{App:matlabCodeBuckling} between lines 156-180. However, in order to reduce looping and avoid repeated calculations, the operations are grouped in a slightly different way. Within a single loop, we extract each of the \texttt{nEig} modes, compute contribution $\circled{2}$ (Line 160-161), build the array $\textsf{p}$ only once (Line 163), and use this for computing contribution $\circled{1}$ (Line 164) and the adjoint loads (Line 168). Then, the adjoint problem is solved simulaneously for all the right hand sides and contribution $\circled{3}$ is computed between lines 173-177. Finally, the sensitivities of the BLFs-related quantities $\mu_{i}$ are given by Line 179.

\subsection{Solution of state and adjoint equations}
\label{sSec:stateAdjointSolvers}

If buckling is among the optimization criteria, two linear systems and an eigenvalue problem have to be solved at each re-design step. It is outside the scope of the present work to discuss how to perform this with highly efficient methods, and we leave this extension to the interested users. Here we resort to built-in Matlab functions, whose efficiency is satisfactory for medium-scale 2D problems.

The stiffness matrix $K$ is decomposed by the Cholesky factorization, using the \texttt{decomposition} Matlab function. The factor \texttt{dK} is stored (Line 119) and used for solving both the linear equilibrium equations (Line 120) and the adjoint problem (Line 173). \texttt{dK} is useful also when calling the \texttt{eigs} function for solving the eigenvalue problem \eqref{eq:eigenvalueEquation}. This is done by the following instructions (see lines 151-152)
\begin{lstlisting}[basicstyle=\scriptsize\ttfamily,breaklines=true,numbers=none,frame=single]
 matFun = @(v) dK\(Kgeo(free,free)*v);
 [eivecs,D]=eigs(matFun,length(free),nEig+4,'sa');
\end{lstlisting}
where the map $\mathbf{v} \mapsto K^{-1}G\mathbf{v}$, repeteadly performed by the iterative Krylov method underlying \texttt{eigs} \citep{lehoucq-soresen_96,stewart_02a}, is explicitly given by the function \texttt{matFun}. This spares the re-factorization of $K$ within \texttt{eigs}, and we have observed a cut of about $10$-$15\%$ of the CPU time, compared to the case where \texttt{matFun} is not provided. When the eigenvalues get close, Krylov based method may loose accuracy on the ``last'' computed eigenpairs \citep{book:wendland18}. Therefore, we heuristically set the number of computed eigenpairs to \texttt{nEig+4}, to retain good accuracy on the \texttt{nEig} used for running the optimization. The argument \texttt{'sa'} selects the smallest algebraic (i.e. ``most negative'') eigenvalues \texttt{D(i)}. Thus, \texttt{-D(i)} correspond to the $\mu_{i}$ of our formulation and, since $\mu_{i} = 1/\lambda_{i}$, to the smallest positive BLFs. Within our assumption of having a load with \textit{fixed} direction it is physically meaningful to accout for positive eigenvalues only. However, the code of \autoref{App:matlabCodeBuckling} could be extended to handle more general loading conditions by properly redefining the eigenvalue equation and the combination of the sensitivity terms in \eqref{eq:sensitivityMu} (see e.g., \cite{lund_09a}).

Finally, a difference with respect to \texttt{top99neo} is that we need to recover the full $K$ and $G$ matrices (lines 118, 149) in order to get correct eigenpairs from \texttt{eigs}. However, the cost of this operation is minor, and essentially identical to the symmetrization operation in \texttt{top88} \citep{andreassen-etal_11a}.

\begin{figure*}
 \centering
  \subfloat[]{
   \includegraphics[scale = 0.25, keepaspectratio]
    {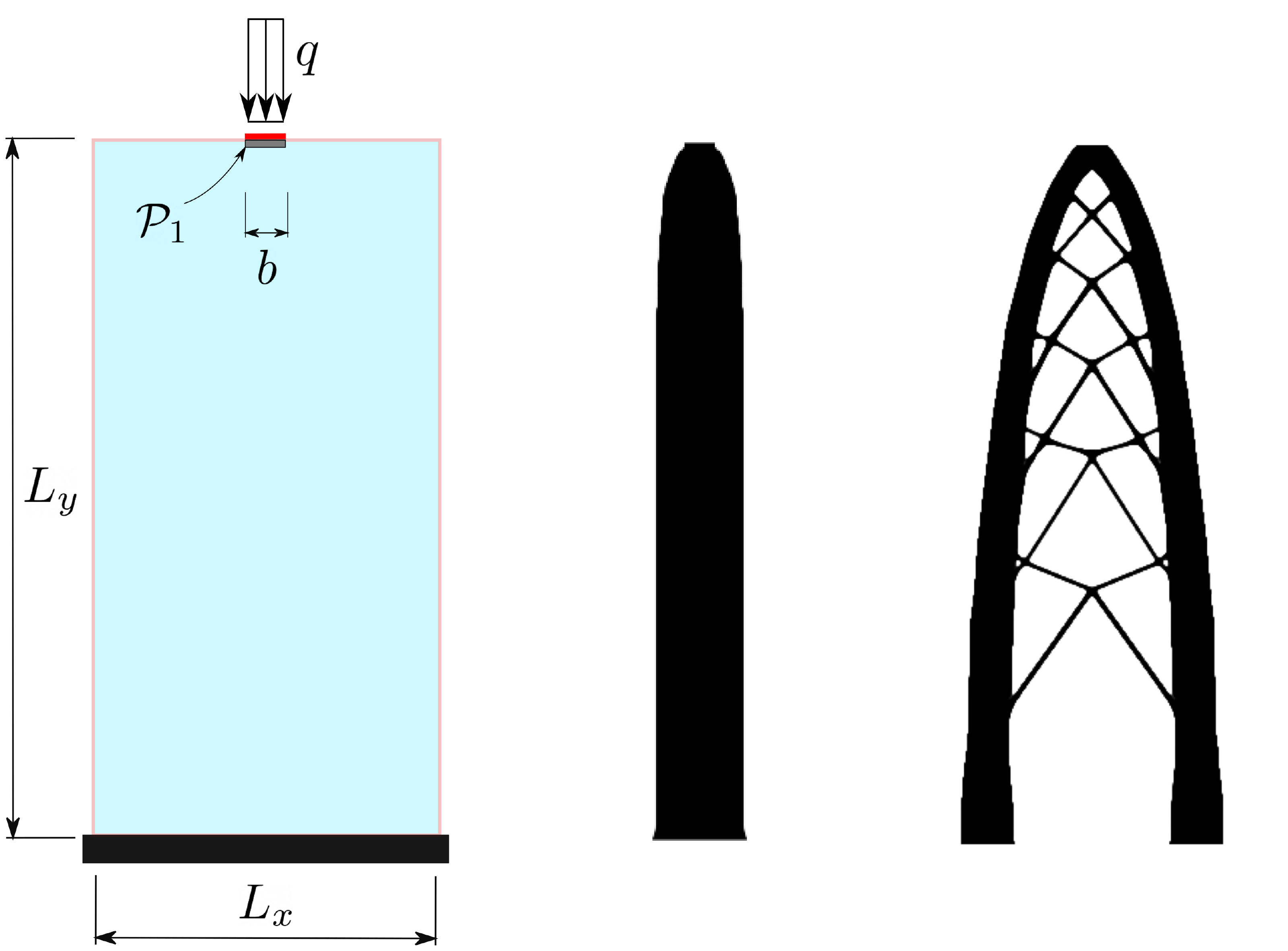}} \qquad
  \subfloat[]{
   \includegraphics[scale = 0.475, keepaspectratio]
    {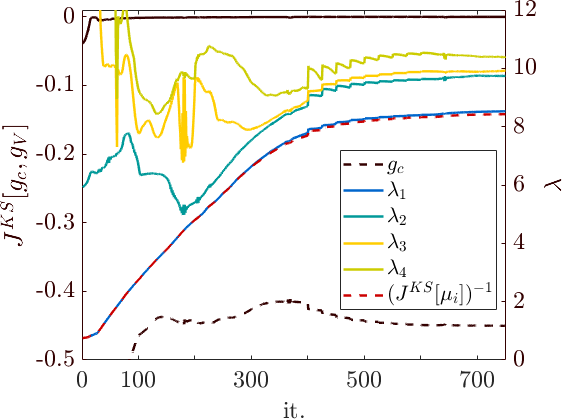}}
 \caption{\small{(a) Geometrical sketch for the compressed column example (leftmost), minimum volume design (center) and buckling reinforced design (rightmost). The minimum volume design has a volume fraction $f \approx 0.24$, compliance $c^{\ast} = 8.57 \cdot 10^{-6}$ and BLF $\lambda_{1} = 0.75$. The buckling design has volume fraction $f = 0.25$ (the volume constraint is tight), compliance $\approx 1.4 c^{\ast}$ and $\lambda_{1} = 8.53$. (b) shows the evolution of the aggregated constraint (essentially coincident with $g_{V}$) and of the compliance constraint $g_{c}$, plotted against the left axis, and of the lowest four BLFs and the reciprocal of the KS function, plotted against the right axis}}
 \label{fig:sketchExamplesCC}
\end{figure*}

\section{Examples}
 \label{Sec:examples}

We present two examples, demonstrating the capabilities of the code provided in \autoref{App:matlabCodeBuckling}. We set $E_{0} = 1$, $E_{\rm min} = 10^{-6}$, $\nu = 0.3$ and we adopt the interpolations \eqref{eq:stiffnessAndStressInterpolations} for all the following problems. The tests have been run on a laptop equipped with an Intel(R) Core(TM) i7-5500U@2.40GHz CPU, 15GB of RAM and Matlab 2018b running in serial mode under Ubuntu 18.04.

\begin{figure}
 \centering
  \includegraphics[scale = 0.4, keepaspectratio]
    {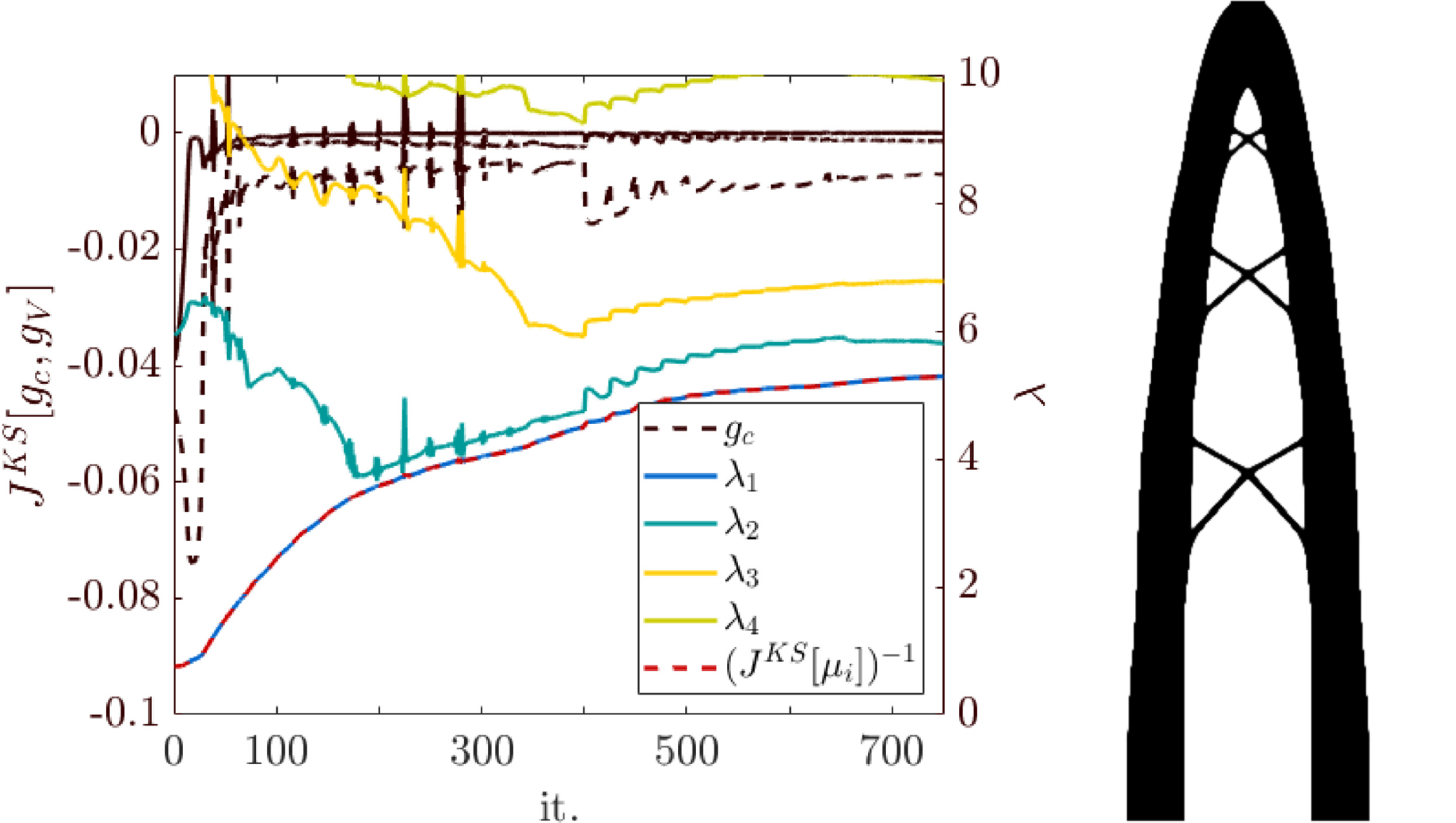}
 \caption{\small{Buckling-reinforced design corresponding to the compliance upper bound $\bar{c}=1.05 c^{\ast}$ and volume fraction $\bar{f}=0.25$ (see \autoref{fig:sketchExamplesCC} for comparison). The final value attained by the fundamental BLF is $\lambda_{1} = 5.28$. The compliance constraint (black dashed curve) now becomes closer to the volume one (black dash-dotted curve), and contributes up to $40\%$ to the KS aggregated constrain value (black continuous curve)}}
 \label{fig:tmpTighter}
\end{figure}

\subsection{BLF maximization of a compressed column}
 \label{sSec:ssBeam}

We consider the cantilever of \autoref{fig:sketchExamplesCC} (a) discretized by $\Omega_{h} = 480 \times 240$ elements. The load is spread over the length $b/L_{y} = 1/15$ at the top and has total magnitude $|q| = 1 \cdot 10^{-3}$. The rectangular area of $10\times 20$ elements near the load constitutes $\mathcal{P}_{1}$ while $\mathcal{P}_{0} = \emptyset$. The boundary conditions and passive regions are implemented by default in the code at lines 64-69, and we remark that the uniform force density $|F_{i}| = 0.625$ is applied on the interior nodes of $b$, whereas $|F_{i}| = 0.3125$ for the two border nodes, for consistency.

We obtain the minimum volume design starting from the fully solid domain ($\mathbf{x}=\mathbf{1}_{m}$) and for the maximum allowed compliance $\bar{c} = 2.5 c_{\mathbf{x} = 1}$, where $c_{\mathbf{x}=1}$ is the initial compliance. Explicitly, we solve
\begin{equation}
 \label{eq:optProblemCantilever_minVol}
 \begin{cases}
 & \min\limits_{\mathbf{x} \in [ 0, 1 ]^{m}}
 f(\mathbf{\hat{x}}) \\
 {\rm s.t.} & g_{c}(\hat{\mathbf{x}}) =
 \mathbf{F}^{T}\mathbf{u}(\hat{\mathbf{x}}) /
 \overline{c} - 1 \leq 0
 \end{cases}
\end{equation}
and we set $p = 3$, fixed throughout the optimization and $r_{\rm min} = 4$ elements for the density filter. The projection parameters are $\eta = 0.5$, $\beta = 2$ and the continuation scheme on $\beta$ follows the rule \texttt{cntBeta=\{150,12,25,2\}}. The code in \autoref{App:matlabCodeBuckling} is called as follows
\begin{lstlisting}[basicstyle=\scriptsize\ttfamily,breaklines=true,numbers=none,frame=single]
topBuck250(480,240,3,4,2,'N',0.5,2,[0.1,0.7,1.2],
300,2,penalG,nEig,pAgg,{['V','C'],2.5})
\end{lstlisting}
where \texttt{penalG,nEig,pAgg} can take any value, because all the buckling-related operations are skipped.

The optimized design, shown in \autoref{fig:sketchExamplesCC} (a), has a volume fraction $f^{\ast}\approx 0.24$ and meets the compliance requirement ($g^{\ast}_{c} \approx -1\cdot 10^{-7}$, thus the constraint is active), but is clearly very prone to lateral buckling. Indeed, we have $\lambda_{1} \approx 0.75$, which means that this design would buckle under the applied load magnitude.

Therefore, we reinforce the design by solving the BLF maximization problem \eqref{eq:optProblemMaxBLF} with volume fraction and compliance upper bounds $\bar{f} = 0.25$ and $\bar{c} = 2.5 c^{\ast}$, respectively, where $c^{\ast} = 8.57 \cdot 10^{-6}$ is the compliance attained by the minimum volume design. The lowest 12 BLFs are aggregated in the KS function, and the aggregation parameter is set to $\rho = 160$, constant throughout the process. Such a high value of the aggregation parameter is no harm to the numerical stability of the KS function \eqref{eq:ksAggregationLambda1} \citep{raspanti-etal_00a}, and is needed to keep a good approximation, especially for the aggregation of constaints. The penalization and filter parameters are choosen as before, but now we adopt $\beta = 6$, and the continuation on this parameter is changed to \texttt{betaCnt=\{400,24,25,2\}}. In this context, a higher $\beta$ parameter was chosen as a simple way to attenuate the artificial buckling mode phenomenon \citep{neves-etal_02a}. Further analyses on this topic, and the implementation of more advanced remedies are left to the interested users.

Assuming that the minimum volume design has been saved in the file \texttt{'IG.mat'}, we call
\begin{lstlisting}[basicstyle=\scriptsize\ttfamily,breaklines=true,numbers=none,frame=single]
topBuck250(480,240,3,4,2,'N',0.5,2,[0.1,0.7,1.2],
750,2,3,12,200,{['B','C','V'],[2.5,0.25]})
\end{lstlisting}

The final design, shown in \autoref{fig:sketchExamplesCC} (a), carries over 10 times more load before buckling than the initial one ($\lambda_{1}$ increases from 0.75 to 8.53), and this comes at a slight increase of volume ($\approx 4\%$), and a stiffness reduction of about $30\%$; however, the compliance constraint never becomes active. \autoref{fig:sketchExamplesCC} (b) shows the evolution of the aggregated constraint, of the (inverse) objective and of the lowest 4 BLFs. Due to the rather high $\rho$ value, the lower bound given by $(J^{KS}_{0})^{-1}$ is a good approximation to $\lambda_{1}$ in the whole optimization history.

\autoref{fig:tmpTighter} shows the design obtained when reducing the maximum allowed compliance to $\bar{c} = 1.05c^{\ast}$. The compliance constraint now becomes much more important, and contributes $30$-$45\%$ to the overall constraint sensitivity, thorughout the optimization. Because of the tightened compliance constraint allowing less freedom \citep{gao-ma_15a,ferrari-sigmund_19a}, the fundamental BLF now attains only a lower value of $\lambda_{1} = 5.28$. 

\autoref{fig:tmpTighter} and other tests we have performed indicate that the simple OC scheme implemented in \texttt{ocUpdate} can cope with the situation where both constraints become very close. However, oscillations in the convergence history may happen as the two constraints become (almost) active, and this may hamper the robustness of the update scheme, especially if the volume-preserving projection is selected (\texttt{ft=3}), and high $\beta$ values are used. In such cases, users are encouraged to shift to more robust optimizers, such as the MMA (see \autoref{App:ksProperties&OCupdate}).

Concerning the computational cost, the whole optimization process (750 re-design steps) takes about $12,900 s$ ($\approx 17.5 s$ per iteration). Of this time, $85.4\%$ ($\approx 11,000 s$) is spent solving the eigenvalue problem, $5.8\%$ ($\approx 750 s$) solving the linear equilibrium equations and $1.1\%$ ($\approx 146.5 s$) solving the adjoint problem. All the other operations are very cheap and, in particular, the overall time fraction spent on the construction of the stress stiffness matrix is $0.3\%$ ($\approx 48 s$) and that for the sensitivity analysis is $5\%$ ($\approx 665.2 s$). The overall time spent on the \texttt{ocUpdate} (see \autoref{App:ksProperties&OCupdate}) for updating the design variables is $0.6\%$ ($\approx 78 s$).

The time spent on the state and adjoint analyses, representing the vast majority of the overall one, can be largely cut by using efficient multi-level solvers for both linear and eigenvalue equations, and this would make buckling optimization with the \texttt{topBuck250} essentially as efficient as multiple-load compliance optimization \citep{ferrari-sigmund_20a}.

\subsection{Wall-reinforcement problem}
 \label{sSec:wallReinforcement}

\begin{figure}
 \centering
  \includegraphics[scale = 0.275, keepaspectratio]
    {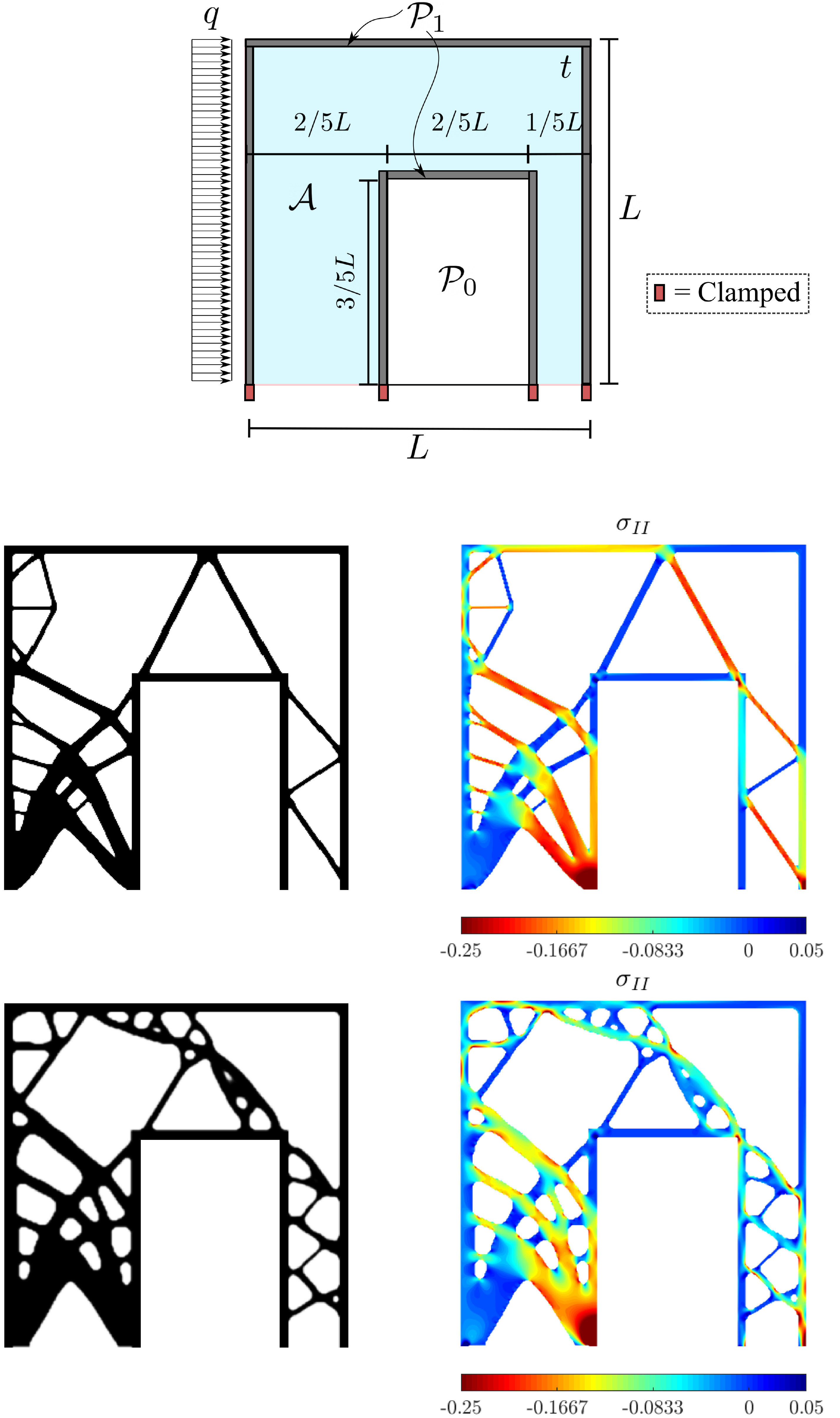}
 \caption{\small{Geometrical setup for the wall reinforcement problem (top), minimum volume design obtained with the compliance constraint alone (middle row) and with a lower bound $\underline{\lambda} = 1.05$ on the BLF (bottom row). Plots on the right column show the distribution of the minimum principal stresses $\sigma_{II}$}}
 \label{fig:sketchExamplesWR}
\end{figure}

\autoref{fig:sketchExamplesWR} shows a wall with an opening (represented by the passive void region $\mathcal{P}_{0}$). The outer frame of the wall and of the opening are both surrounded by a solid frame ($\mathcal{P}_{1}$) with thickness $t = L/40$, clamped at the base. The outer frame alone represents a suitable load transfer path for the applied load $q$, but would result in a structure that is too compliant.

Therefore, the goal is to find the minimum weight configuration of reinforcement material in the active domain $\mathcal{A}$, starting from the full solid design (i.e. $\mathbf{x} = 1$ on $\mathcal{A}$, corresponding to a volume fraction of $\approx 0.76$ of the overall domain $\mathcal{A}\cup\mathcal{P}_{1}\cup\mathcal{P}_{0}$). To do this, we refer to problem \eqref{eq:optProblemMinVolcstrBLFandC}, where $\bar{c}$ is again set to $2.5$ times the compliance of the full design, and the BLF is kept above the value $\underline{\lambda}$, by applying the constraints $g_{c}$ and $g_{\lambda}$ of \eqref{eq:constraintsForm}.

We set $L=1$, $\Omega_{h} = 320\times 320$ and the load, uniformly distributed over the whole left edge, has total magnitude $|q| = 10^{-2}$. To implement this example we only need to replace lines 64-70 with the following
\begin{lstlisting}[basicstyle=\tiny\ttfamily,breaklines=true,numbers=none,frame=single]
fixedNodes=nodeNrs(end,[1:9,120:128,257:265,313:321])
fixed=[2*fixedNodes,2*fixedNodes-1];
lcDof=2*nodeNrs(1:end,1)-1;
modF=1e-2/Ly/(length(lcDof)-1);
F=fsparse(lcDof,1,modF,[nDof,1]);
[F(lcDof(1)),F(lcDof(end))]=deal(F(lcDof(1))/2,F(lcDof(end))/2);
a1=elNrs(1:nely/40,:);
a2=elNrs(:,[1:nelx/40,end-nelx/40+1:end]);
a3=elNrs(2*nely/5:end,2*nelx/5:end-nelx/5);
a4=elNrs(2*nely/5-nely/40:end,2*nelx/5-nelx/40:2*nelx/5-1);
a5=elNrs(2*nely/5-nely/40:end,end-nelx/5+1:end-nelx/5+nelx/40);
a6=elNrs(2*nely/5-nely/40:2*nely/5-1,2*nelx/5:end-nelx/5);
[pasS,pasV]=deal(unique([a1(:);a2(:);a4(:);a5(:);a6(:)]), a3(:));
\end{lstlisting}

The other parameters for the optimization are chosen as follows: $p_{K}$ and $p_{G}$ are increased from 3 to 6, both according to the continuation scheme \texttt{\{25,6,25,0.25\}}, $\eta = 0.5$ and $\beta$ is increased from 2 to 12 according to the continuation scheme \texttt{betaIncrease=\{325,12,25,2\}}. The density filter radius is $r_{\rm min} = 3$ and the lowest 12 BLFs are aggregated in the KS function, with $\rho = 160$.

The middle row of \autoref{fig:sketchExamplesWR} shows the optimized design obtained by calling the code in \autoref{App:matlabCodeBuckling} as
\begin{lstlisting}[basicstyle=\scriptsize\ttfamily,breaklines=true,numbers=none,frame=single]
topBuck250(320,320,3,3,2,'N',0.5,2,[0.1,0.7,1.2],
500,1,penalG,nEig,pAgg,{['V','C'],2.5})
\end{lstlisting}
thus, without the buckling constraint ($\underline{\lambda} = 0$). For this design, the compliance constraint is tight ($g^{\ast}_{c} \approx 10^{-9}$) and the volume fraction attains the value $f^{\ast} = 0.245$.

\begin{figure}[tb]
 \centering
  \includegraphics[scale = 0.315, keepaspectratio]
    {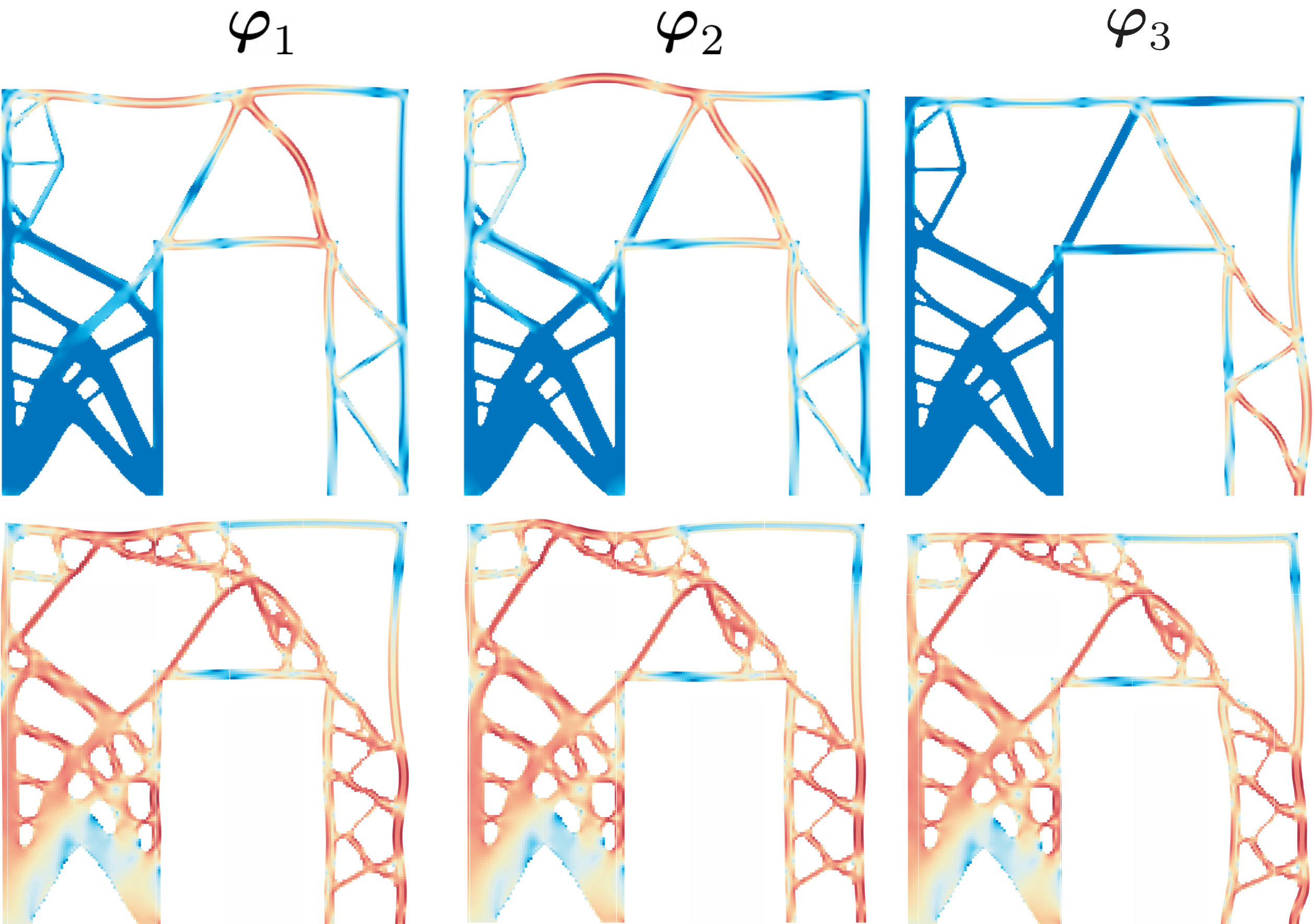}
 \caption{\small{The lowest three buckling modes ($\boldsymbol{\varphi}_{i}$) for the wall reinforced design, obtained without (top row) and with (bottom row) the buckling constraint. The colorscale (blue to red) follows the logarithm of the strain energy density (low to high)}}
 \label{fig:bucklingModesWR}
\end{figure}

The reinforcement material is concentrated in the lower left part of the structure, assuming a configuration that vaguely resembles that of an optimized cantilever for a tip load. Elsewhere, the reinforcement material is organized in slender bars. This design is clearly not effective against buckling, as can be seen from the distribution of the minimum principal stresses ($\sigma_{II}$) in the rightmost plot of \autoref{fig:sketchExamplesWR}. From this plot we acknowledge that some thin bars are subject to relatively high compression, and from \autoref{fig:bucklingModesWR} we see that they are indeed undergoing localized buckling. The fundamental BLF is $\lambda_{1} = 3.15$ for the initial, solid structure and $\lambda_{1} = 0.15$ for the minimum volume design. Therefore, this latter is of little practical use.

To obtain a more realistic design, we solve again the optimization problem imposing the lower bound $\underline{\lambda} = 1.05$ to the BLF, calling the routine by
\begin{lstlisting}[basicstyle=\scriptsize\ttfamily,breaklines=true,numbers=none,frame=single]
topBuck250(320,320,3,3,2,'N',0.5,2,[0.1,0.7,1.2],
500,1,3,12,160,{['V','C','B'],[2.5,1.05]})
\end{lstlisting}
and the continuation on $p_{G}$ is identical to that on $p_{K}$.

We obtain the design shown in the bottom row of \autoref{fig:sketchExamplesWR}, having a volume fraction $f^{\ast} = 0.33$ and fulfilling both constraints. The minimum principal stress ($\sigma_{II}$) plot shows that the thinnest bars are now subjected to tension only, whereas compressed parts have a more organic, frame-like configuration. By comparing the distribution of the compressive stresses with the lowest three buckling modes depicted in \autoref{fig:bucklingModesWR}, we observe that now buckling is triggered by a deformation in the outer frame, and this endows the structure with a much better overall behavior.

\section{Conclusions}
 \label{Sec:conclusions}

We have presented a compact Matlab implementation that allows for solution of buckling TO problems on a laptop, up to a scale that is limited only by the solution of the eigenvalue problem, which is now absorbing the vast majority of the CPU time. The setup of the stress stiffness matrix and the sensitivity analysis of the buckling load factor are made very cheap by rethinking their standard construction and using a fully vectorized implementation. Otherwise, these operations alone would make solution of the problem in Matlab very slow.

The code provided in \autoref{App:matlabCodeBuckling} can be immediately used to solve problems for the maximization of the buckling load factor with compliance and volume constraints, and for the minimization of volume with buckling and compliance constraints. The basic TO problems of volume constrained compliance minimization and \emph{vice versa} can also be run avoiding all the buckling related operations; thus obtaining a code with the same efficiency as \texttt{top99neo}. However, the code can easily be extended to more general problems involving these three criteria, and to implement more advance formulations and solution methods. As an example, methods to deal with the artificial modes phenomenon, here overlooked, can easily be implemented. Also, preconditioned eigensolvers \citep{dunning-etal_16a} and multilevel methods \citep{ferrari-sigmund_20a} are easy to implement in the present code, to alleviate the burden of the eigenvalue solution. Finally, we point out that other problems that involve element-wise operations, such as non-linear elasticity or orientation dependent material properties, may be effectively reformulated with a strategy similar to what is described here.

Therefore, we believe that the present code can be a helpful contribution to researchers and practitioners who aim at solving TO problems with buckling criteria in Matlab, and provide food for thought to anyone seeking efficient implementations of potentially time consuming problems in FE analysis and TO.

\begin{acknowledgements}
 The authors FF and JKG acknowledge the support of the National Aeronautics and Space Administration (NASA) under Grant No. 80NSSC18K0428. The opinions, findings, and conclusions or recommendations expressed in this paper are those of the authors and do not necessarily reflect the views of NASA. OS is supported by the Villum Fonden through the Villum Investigator Project ``InnoTop''.
\end{acknowledgements}

\begin{small}
\paragraph{\textsc{\textbf{Reproducibility of results}}}
Matlab code is listed in the Appendix and available at \texttt{www.topopt.dtu.dk}. The \texttt{stenglib} package, containing the \texttt{fsparse} function, is avaialble for download at \texttt{https://github.com/stefanengblom/stenglib}.
\end{small}


\begin{small}
\bibliographystyle{spbasic}      
\bibliography{databaseSMO_edubuckling.bib}

\begin{thebibliography}{40}
\providecommand{\natexlab}[1]{#1}
\providecommand{\url}[1]{{#1}}
\providecommand{\urlprefix}{URL }
\expandafter\ifx\csname urlstyle\endcsname\relax
  \providecommand{\doi}[1]{DOI~\discretionary{}{}{}#1}\else
  \providecommand{\doi}{DOI~\discretionary{}{}{}\begingroup
  \urlstyle{rm}\Url}\fi
\providecommand{\eprint}[2][]{\url{#2}}

\bibitem[{Andreassen et~al(2011)Andreassen, Clausen, Schevenels, Lazarov, and
  Sigmund}]{andreassen-etal_11a}
Andreassen E, Clausen A, Schevenels M, Lazarov BS, Sigmund O (2011) Efficient
  topology optimization in matlab using 88 lines of code. Structural and
  Multidisciplinary Optimization 43(1):1--16, \doi{10.1007/s00158-010-0594-7}

\bibitem[{Bathe(1982)}]{book:bathe_fem}
Bathe KJ (1982) {F}inite {E}lement procedures in engineering analysis, 1st edn

\bibitem[{Bends\o{}e and Sigmund(2004)}]{book:bendsoe-sigmund_2004}
Bends\o{}e MP, Sigmund O (2004) Topology Optimization: Theory, Methods and
  Applications. Springer

\bibitem[{de~Borst et~al(2012)de~Borst, Crisfield, Remmers, and
  Verhoosel}]{book:deborst2012}
de~Borst R, Crisfield MA, Remmers JJC, Verhoosel CV (2012) Non--Linear Finite
  Element Analysis of Solids and Structures, second edition edn. John Wiley \&
  Sons

\bibitem[{Bourdin(2001)}]{bourdin_01a}
Bourdin B (2001) Filters in topology optimization. International Journal for
  Numerical Methods in Engineering 50(9):2143--2158, \doi{10.1002/nme.116}

\bibitem[{Brent(1973)}]{book:Brent73}
Brent R (1973) Algorithms for Minimization without Derivatives. Prentice Hall

\bibitem[{Bruns and Tortorelli(2001)}]{bruns-tortorelli_01a}
Bruns TE, Tortorelli DA (2001) Topology optimization of non-linear elastic
  structures and compliant mechanisms. Computer Methods in Applied Mechanics
  and Engineering 190(26):3443 -- 3459,
  \doi{http://dx.doi.org/10.1016/S0045-7825(00)00278-4},
  \urlprefix\url{http://www.sciencedirect.com/science/article/pii/S0045782500002784}

\bibitem[{Bruyneel et~al(2008)Bruyneel, Colson, and
  Remouchamps}]{bruyneel-etal_08a}
Bruyneel M, Colson B, Remouchamps A (2008) Discussion on some convergence
  problems in buckling optimisation. Structural and Multidisciplinary
  Optimization 35(2):181--186

\bibitem[{Crisfield(1991)}]{book:crisfield91}
Crisfield MA (1991) Nonlinear finite element analysis of solid and structures,
  vol~I. Wiley

\bibitem[{Dunning et~al(2016)Dunning, Ovtchinnikov, Scott, and
  Kim}]{dunning-etal_16a}
Dunning PD, Ovtchinnikov E, Scott J, Kim A (2016) Level--set topology
  optimization with many linear buckling constraints using and efficient and
  robust eigensolver. International Journal for Numerical Methods in
  Engineering

\bibitem[{Engblom and Lukarski(2016)}]{engblom-lukarski_16a}
Engblom S, Lukarski D (2016) Fast matlab compatible sparse assembly on
  multicore computers. Parallel Computing 56:1--17

\bibitem[{Ferrari and Sigmund(2019)}]{ferrari-sigmund_19a}
Ferrari F, Sigmund O (2019) Revisiting topology optimization with buckling
  constraints. Structural and Multidisciplinary Optimization 59(5):1401--1415,
  \doi{10.1007/s00158-019-02253-3}

\bibitem[{Ferrari and Sigmund(2020{\natexlab{a}})}]{ferrari-sigmund_20b}
Ferrari F, Sigmund O (2020{\natexlab{a}}) A new generation 99 line matlab code
  for compliance topology optimization and its extension to 3d. Structural and
  Multidisciplinary Optimization 62:2211--2228

\bibitem[{Ferrari and Sigmund(2020{\natexlab{b}})}]{ferrari-sigmund_20a}
Ferrari F, Sigmund O (2020{\natexlab{b}}) Towards solving large-scale topology
  optimization problems with buckling constraints at the cost of linear
  analyses. Computer Methods in Applied Mechanics and Engineering 363:112,911,
  \doi{https://doi.org/10.1016/j.cma.2020.112911}

\bibitem[{Gao and Ma(2015)}]{gao-ma_15a}
Gao X, Ma H (2015) Topology optimization of continuum structures under buckling
  constraints. Computers \& Structures 157:142--152

\bibitem[{Gravesen et~al(2011)Gravesen, Evgrafov, and
  Nguyen}]{gravesen-etal_11a}
Gravesen J, Evgrafov A, Nguyen DM (2011) On the sensitivities of multiple
  eigenvalues. Structural and Multidisciplinary Optimization 44(4):583--587

\bibitem[{Groenwold and Etman(2008)}]{groenwold-etman_08a}
Groenwold AA, Etman LFP (2008) On the equivalence of optimality criterion and
  sequential approximate optimization methods in the classical topology layout
  problem. Internal Journal for Numerical Methods in Engineering 73:297--316

\bibitem[{Guest et~al(2004)Guest, Pr\'{e}vost, and Belytschko}]{guest-etal_04a}
Guest JK, Pr\'{e}vost JH, Belytschko T (2004) Achieving minimum length scale in
  topology optimization using nodal design variables and projection functions.
  International Journal for Numerical Methods in Engineering 61(2):238--254,
  \doi{10.1002/nme.1064}

\bibitem[{Guest et~al(2011)Guest, Asadpoure, and a}]{guest-etal_11a}
Guest JK, Asadpoure A, a SH (2011) Eliminating beta continuation from heaviside
  projection and density filter algorithms. Structural and Multidisciplinary
  Optimization 44(4):443--453

\bibitem[{Kennedy and Hicken(2015)}]{kennedy-hicken_15a}
Kennedy GJ, Hicken JE (2015) Improved constraint-aggregation methods. Computer
  Methods in Applied Mechanics and Engineering 289(Supplement C):332 -- 354

\bibitem[{Kreisselmeier and Steinhauser(1979)}]{kreisselmeier-steinhauser_79a}
Kreisselmeier G, Steinhauser R (1979) Systematic control design by optimizing a
  vector performance index. IFAC Proceedings Volumes 12(7):113 -- 117, iFAC
  Symposium on computer Aided Design of Control Systems, Zurich, Switzerland,
  29-31 August

\bibitem[{Lehoucq and Sorensen(1996)}]{lehoucq-soresen_96}
Lehoucq R, Sorensen DC (1996) Deflation techniques for an implicitly re-started
  arnoldi iteration. SIAM J Matrix Anal Appl 17:789--821

\bibitem[{Lund(2009)}]{lund_09a}
Lund E (2009) Buckling topology optimization of laminated multi--material
  composite shell structures. Composite Structures 91(2):158--167

\bibitem[{Neves et~al(2002)Neves, Sigmund, and Bends\o{}e}]{neves-etal_02a}
Neves MM, Sigmund O, Bends\o{}e MP (2002) Topology optimization of periodic
  microstructures with a penalization of highly localized buckling modes.
  International Journal for Numerical Methods in Engineering 54(6):809--834

\bibitem[{Raspanti et~al(2000)Raspanti, Bandoni, and
  Biegler}]{raspanti-etal_00a}
Raspanti CG, Bandoni JA, Biegler LT (2000) New strategies for flexibility
  analysis and desing under uncertainties. Computers \& Chemical Engineering
  24:2193--2209

\bibitem[{Rodrigues et~al(1995)Rodrigues, Guedes, and
  Bends\o{}e}]{rodrigues-etal_95a}
Rodrigues HC, Guedes JM, Bends\o{}e MP (1995) Necessary conditions for optimal
  design of structures with a nonsmooth eigenvalue based criterion. Structural
  Optimization 9:52--56

\bibitem[{Seyranian et~al(1994)Seyranian, Lund, and Olhoff}]{seyranian_94a}
Seyranian AP, Lund E, Olhoff N (1994) Multiple eigenvalues in structural
  optimization problems. Structural optimization 8(4):207--227

\bibitem[{Sigmund(2001)}]{sigmund_01a}
Sigmund O (2001) A 99 line topology optimization code written in {M}atlab.
  Structural and Multidisciplinary Optimization 21(2):120--127,
  \doi{10.1007/s001580050176}

\bibitem[{Sigmund(2007)}]{sigmund_07a}
Sigmund O (2007) Morphology--based black and white filters for topology
  optimization. Structural and Multidisciplinary Optimization 33(4):401--424

\bibitem[{Stewart(2002)}]{stewart_02a}
Stewart G (2002) A {K}rylov-{S}chur algorithm for large eigenproblems. SIAM
  Journal on Matrix Analysis and Applications 23(3):601--614

\bibitem[{Svanberg(1987)}]{svanberg_87a}
Svanberg K (1987) The method of moving asymptotes - {A} new method for
  structural optimization. International Journal for Numerical Methods in
  Engineering 24(2):359--373

\bibitem[{Svanberg(2002)}]{svanberg_02a}
Svanberg K (2002) A class of globally convergent optimization methods based on
  conservative convex separable approximations. SIAM Journal on Optimization
  12(2):555--573, \doi{10.1137/S1052623499362822}

\bibitem[{Svanberg(2007)}]{svanberg-note}
Svanberg K (2007) {MMA} and {GCMMA}--two methods for nonlinear optimization,
  available for download at https://people.kth.se/~krille/mmagcmma.pdf

\bibitem[{Wahlbin(1995)}]{book:wahlbin1995}
Wahlbin LB (1995) Superconvergence in galerkin finite element methods

\bibitem[{Wang et~al(2011)Wang, Lazarov, and Sigmund}]{wang-etal_11a}
Wang F, Lazarov B, Sigmund O (2011) On projection methods, convergence and
  robust formulations in topology optimization. Structural and
  Multidisciplinary Optimization 43(6):767--784

\bibitem[{Wendland(2018)}]{book:wendland18}
Wendland H (2018) Numerical linear algebra: {A}n introduction. Cambridge texts
  in applied mathematics, Cambridge University Press

\bibitem[{Wrenn(1989)}]{wrenn_89a}
Wrenn G (1989) An indirect method for numerical optimization using the
  {K}reisselmeir--{S}teinhauser function. Tech. rep.,
  \urlprefix\url{https://books.google.it/books?id=LjsCAAAAIAAJ}

\bibitem[{Xia and Breitkopf(2015)}]{xia-breitkopf_15a}
Xia L, Breitkopf P (2015) Design of materials using topology optimization and
  energy-based homogenization approach in matlab. Structural and
  Multidisciplinary Optimization 52(6):1229--1241,
  \doi{10.1007/s00158-015-1294-0}

\bibitem[{Yago et~al(2020)Yago, Cante, Lloberas-Valls, and
  Oliver}]{yago-etal_20a}
Yago D, Cante J, Lloberas-Valls O, Oliver J (2020) Topology optimization using
  the unsmooth variational topology optimization (unvartop) method: an
  educational implementation in matlab. Structural and Multidisciplinary
  Optimization Under review

\bibitem[{Zienkiewicz and Taylor(2005)}]{book:zienkiewicz-taylor06}
Zienkiewicz O, Taylor R (2005) The finite element method for solid and
  structural mechanics, 6th edn. Elsevier

\end{thebibliography}
\end{small}

\appendix

\section{Notes}
 \label{App-notesDiscretizationOperators}

Let $X^{e} = [x^{e}_{i},y^{e}_{i}]^{k}_{i=1}$ be the array collecting the coordinates of a $k$-noded element and $N(\xi,\zeta) = [ N_{1}(\xi,\zeta), N_{2}(\xi,\zeta), \ldots, N_{k}(\xi,\zeta)]$ be the array collecting the shape functions, where $(\xi,\zeta)$ are the logical coordinates \citep{book:bathe_fem}. The $(x,y)$-gradient of the shape functions is computed as $\nabla_{(x, y)}N = J^{-1} \nabla_{(\xi, \zeta)}N$, where the $(\xi,\zeta)$-gradient and Jacobian are
\begin{equation}
 \label{eq:dNandJ}
 \begin{aligned}
 \nabla_{(\xi, \zeta)}N & = \left[
 \begin{array}{ccc}
 \partial_{\xi} N_{1}, \ldots, \partial_{\xi} N_{k} \\
 \partial_{\zeta} N_{1}, \ldots,
 \partial_{\zeta} N_{k}
 \end{array}
 \right]
 \\
 J(\xi, \eta) & = \left[
 \begin{array}{cc}
 \sum_{k}\partial_{\xi} N_{k} x_{k} &
 \sum_{k}\partial_{\xi} N_{k} y_{k} \\
 \sum_{k}\partial_{\zeta} N_{k} x_{k} &
 \sum_{k}\partial_{\zeta} N_{k} y_{k}
 \end{array}
 \right] = \nabla_{(\xi, \zeta)}N X^{e}
 \end{aligned}
\end{equation}

The discretized deformation gradient ($B_{1}$) and strain-displacement operator ($B_{0}$) are then recovered as
\begin{equation}
 \label{eq:B0B1definition}
 \begin{aligned}
  B_{1} & = \left[
  \begin{array}{c}
  \nabla_{(x, y)}N \otimes [ 1, 0 ]^{T} \\
  \nabla_{(x, y)}N \otimes [ 0, 1 ]^{T}
  \end{array}
  \right]
  \\
  B_{0} & = L ( \nabla_{(x, y)}N \otimes I_{2} )
 \end{aligned}
\end{equation}
where $L$ is the placement matrix (here shown for a $\mathcal{Q}_{4}$ element)
\[
 L = \left[
 \begin{array}{cccc}
 1 & 0 & 0 & 0 \\
 0 & 0 & 0 & 1 \\
 0 & 1 & 1 & 0
 \end{array}
 \right]
\]
and $I_{2}$ is the identity matrix of order 2.

The operations above have a remarkably compact implementation in the code of \autoref{App:matlabCodeBuckling}, for the particular case of a $\mathcal{Q}_{4}$ element (see lines extracted below).
\begin{lstlisting}[basicstyle=\scriptsize\ttfamily,breaklines=true,numbers=none,frame=single]
 lMat=zeros(3,4);lMat(1,1)=1;lMat(2,4)=1;lMat(3,2:3)=1;
 dN=@(xi,zi)0.25*[zi-1,1-zi,1+zi,-1-zi;xi-1,-1-xi,1+xi,1-xi];
 gradN=(dN(xi,zi)*xe)/dN(xi,zi)                               
 B0=@(gradN) lMat*kron(gradN,eye(2))
 B1=[kron(gradN,[1,0]);kron(gradN,[0,1])]
\end{lstlisting}

The instructions above could be further simplified, since we are considering a uniform discretization (i.e. $J$ is constant). However, the current definition is general enough to be easily extended to higher order elements by only modifying the definition of \texttt{dN}.

\onecolumn

\section{Details on the re-design rule and on the KS aggregation function}
 \label{App:ksProperties&OCupdate}

Let us have the two sets of functions $\{g_{0(i)}(\mathbf{x})\}_{i=1\ldots,q}$ and $\{g_{j}(\mathbf{x})\}_{j=1\ldots,s}$ and the following optimization problem
\begin{equation}
 \begin{cases}
   & \min\limits_{\mathbf{x}\in\left[ 0, 1 \right]^{m}}
   \max\limits_{i=1,\ldots, q} g_{0(i)}(\mathbf{x}) \\
   {\rm s.t.} & g_{j}\left( \mathbf{x} \right) \leq 0
   \:, \qquad\qquad j = 1, \ldots, s
  \end{cases}
\end{equation}

We replace $J^{KS}_{0}\left( \mathbf{x} \right) := J^{KS}[g_{0(i)}\left( \mathbf{x} \right)]$ to the ``$\max$'' operator in the objective, and $J^{KS}_{1}\left( \mathbf{x} \right) := J^{KS}[g_{j}\left( \mathbf{x} \right)]$ to the set of constraints, obtaining the single objective, single constraint problem
\begin{equation}
 \label{eq:App-optProblemWithKS}
 \begin{cases}
  & \min\limits_{\mathbf{x}\in\left[ 0, 1 \right]^{m}}
  J^{KS}_{0}\left( \mathbf{x} \right) \\
  {\rm s.t.} & J^{KS}_{1}\left( \mathbf{x} \right) \leq 0
 \end{cases}
\end{equation}
which is analogous to \eqref{eq:optProblemMaxBLF}. At the current design point ($\mathbf{x}_{k} = \boldsymbol{\xi}$) the objective and constraint functions are expanded in terms of intervening variables $y_{e}(x_{e})$, $e=1, \ldots m$, each selected such to build a convex approximation. Choosing an MMA-like form for the intervening variables \citep{svanberg_87a}, the objective ($i=0$) and constraint ($i=1$) are expanded as
\begin{equation}
 \label{eq:mmaLikeExpansion}
 \begin{aligned}
 J^{KS}_{i}(\mathbf{x}) & \approx
 J^{KS}_{i}(\boldsymbol{\xi}) +
 \sum^{m}_{e=1} \left[ (U_{e}-\xi_{e})^{2}p^{i}_{e}(\boldsymbol{\xi}) \left( 
 \frac{1}{U_{e}-x_{e}} - \frac{1}{U_{e}-\xi_{e}} \right) - 
 (\xi_{e}-L_{e})^{2}q^{i}_{e}(\boldsymbol{\xi}) \left( 
 \frac{1}{x_{e}-L_{e}} - \frac{1}{\xi_{e}-L_{e}} \right) \right] \\
 & \underbrace{J^{KS}_{i}(\boldsymbol{\xi}) - \sum^{m}_{e=1} \left[
 (U_{e}-\xi_{e}) p^{i}_{e}(\boldsymbol{\xi}) - (\xi_{e}-L_{e}) q^{i}_{e}(\boldsymbol{\xi})\right]}_{\widehat{J^{KS}_{i}}(\boldsymbol{\xi})} +
 \sum^{m}_{e=1} \left[ \frac{(U_{e}-\xi_{e})^{2}p^{i}_{e}(\boldsymbol{\xi})}
 {U_{e}-x_{e}} -
 \frac{(\xi_{e}-L_{e})^{2}q^{i}_{e}(\boldsymbol{\xi})}{x_{e}-L_{e}} \right] \\
 & = \widehat{J^{KS}_{i}}(\boldsymbol{\xi}) + 
 \sum^{m}_{e=1} \left[\frac{(U_{e}-\xi_{e})^{2}
 p^{i}_{e}(\boldsymbol{\xi})}{U_{e}-x_{e}} -
 \frac{(\xi_{e}-L_{e})^{2}
 q^{i}_{e}(\boldsymbol{\xi})}{x_{e}-L_{e}} \right]
 \end{aligned}
\end{equation}

where $p^{i}_{e}(\boldsymbol{\xi}) = \max\{\partial_{e}J^{KS}_{i}(\boldsymbol{\xi}),0\} \geq 0$ and $q^{i}_{e}(\boldsymbol{\xi}) = \min\{\partial_{e}J^{KS}_{i}(\boldsymbol{\xi}),0\} \leq 0$ and $(L_{e}, U_{e})$, $e=1,\ldots, m$ are the moving asymptotes, such that $L_{e} < (x_{e}, \xi_{e}) < U_{e}$.

Once \eqref{eq:mmaLikeExpansion} is substituted into \eqref{eq:App-optProblemWithKS}, the Lagrangian of the problem reads
\begin{equation}
 \label{eq:lagrangianFull}
 \ell( \mathbf{x}, \lambda ) = 
 \widehat{J^{KS}_{0}}(\boldsymbol{\xi}) +
 \kappa
 \widehat{J^{KS}_{1}}(\boldsymbol{\xi}) +
 \sum^{m}_{e=1} \left[ \frac{(U_{e}-\xi_{e})^{2}
 (p^{0}_{e}(\boldsymbol{\xi})+\kappa
 p^{1}_{e}(\boldsymbol{\xi}))
 }{U_{e}-x_{e}} -
 \frac{(\xi_{e}-L_{e})^{2}
 (q^{0}_{e}(\boldsymbol{\xi})+\kappa
 q^{1}_{e}(\boldsymbol{\xi}))}{x_{e}-L_{e}}
 \right]
\end{equation}
where $\kappa \geq 0$ is the Lagrange multiplier. The stationarity condition of the Lagrangian with respect to $x_{e}$
\begin{equation}
 \label{eq:stationarityLagrangian}
 \partial_{e} \ell( \mathbf{x}, \lambda ) = 
 \frac{(U_{e}-\xi_{e})^{2}
 (p^{0}_{e}(\boldsymbol{\xi})+\kappa
 p^{1}_{e}(\boldsymbol{\xi}))
 }{(U_{e}-x_{e})^{2}} +
 \frac{(\xi_{e}-L_{e})^{2}
 (q^{0}_{e}(\boldsymbol{\xi})+\kappa
 q^{1}_{e}(\boldsymbol{\xi}))}{(x_{e}-L_{e})^{2}} = 0
\end{equation}
can be solved explicitly, giving the primal update map
\begin{equation}
 \label{eq:projectionFmma}
 x_{e}(\kappa) = 
 \max\left\{\delta_{-},\min\left\{ \delta_{+},
 \frac{L_{e}(U_{e}-\xi_{e})\sqrt{p^{0}_{e}(\boldsymbol{\xi}) + \kappa p^{1}_{e}(\boldsymbol{\xi})} + U_{e}(\xi_{e}-L_{e})\sqrt{-q^{0}_{e}(\boldsymbol{\xi}) - \kappa q^{1}_{e}(\boldsymbol{\xi})}}
 {(U_{e}-\xi_{e})\sqrt{p^{0}_{e}(\boldsymbol{\xi}) + \kappa p^{1}_{e}(\boldsymbol{\xi})} + (\xi_{e}-L_{e})\sqrt{-q^{0}_{e}(\boldsymbol{\xi}) - \kappa q^{1}_{e}(\boldsymbol{\xi})}}
 \right\}\right\}
\end{equation}
where $0 < \delta_{-} < \delta_{+} < 1$ are the adaptive move limits, depending on the asymptotes and we have highlighted the dependence on the dual variable $\kappa$. This latter can be computed from the following equation, obtained by formally substituting \eqref{eq:projectionFmma} into \eqref{eq:lagrangianFull} and imposing stationarity with respect to $\kappa$
\begin{equation}
 \label{eq:stationarityDualFunction}
 \psi( \kappa ) = 
 \widehat{J^{KS}_{1}}(\boldsymbol{\xi}) +
 \sum^{m}_{e=1} \left[ \frac{(U_{e}-\xi_{e})^{2}
 p^{1}_{e}(\boldsymbol{\xi})
 }{U_{e}-x_{e}(\kappa)} -
 \frac{(\xi_{e}-L_{e})^{2}
 q^{1}_{e}(\boldsymbol{\xi})}{x_{e}(\kappa)-L_{e}}
 \right] = 0
\end{equation}

The re-design procedure is implemented in the routine \texttt{ocUpdate}, listed at the bottom of this Section. The input parameters are the current re-design step (\texttt{loop}), the current design point, objective sensitivity, constraint value and sensitivity (\texttt{xT}, \texttt{dg0}, \texttt{g1}, \texttt{dg1}), the design variables at the two previous iterations (\texttt{xOld}, \texttt{xOld1}), and the current asymptotes (\texttt{as}). The input \texttt{ocPar=[move,asReduce,asRelax]} collects the steplenght (\texttt{move}) and the two numbers used for tightening and relaxing the asymptotes according to the smootheness of the optimization history (\texttt{asReduce,asRelax}).

The evolution of the asymptotes  and their link with the adaptive bounds $\delta_{-}$ and $\delta_{+}$ (called \texttt{xL} and \texttt{xU}) are conceptually identical to that suggested by \cite{svanberg_87a} (see lines 4-13). Following \cite{guest-etal_11a} we allow the tightening of the initial asymptotes to accomodate large $\beta$ values, by passing this parameter to the \texttt{ocUpdate} routine. The user can also reset the asymptotes whenever continuation is applied on the penalization and projection parameters, by passing the additional variable \texttt{restartAsy}.

The primal update \eqref{eq:projectionFmma} and stationarity condition \eqref{eq:stationarityDualFunction} are defined on lines 20-23, based on the positive variables \texttt{p0}, \texttt{q0}, \texttt{p1} and \texttt{q1}, representing $(U_{e}-\xi_{e})^{2}p^{i}_{e}$ and $-(\xi_{e}-L_{e})^{2}q^{i}_{e}$, $i=0,1$, respectively (lines 17-18). The dual variable $\kappa$, is computed between lines 26-33 accounting for three possible situations: given the search window $[0, \bar{\kappa}]$ (e.g. $\bar{\kappa} = 10^{6}$)
\begin{enumerate}
 \item If $\psi(0)\psi(\bar{\kappa}) < 0 \rightarrow$ the Lagrange multiplier is within $(0, \bar{\kappa})$, and we compute it by the built-in Matlab function \texttt{fzero}, applying a version of the Brent's algorithm \citep{book:Brent73} (lines 27-28);
 \item $\psi(0)<0 \rightarrow$ the constraint is not active within the current local expansion. Thus, we set $\kappa = 0$ and $x_{e} = x_{e}(0)$ (line 30);
 \item $\psi(\bar{\kappa})>0 \rightarrow$ the constraint cannot be fulfilled within the current local expansion (i.e. we cannot find a feasible solution). In this case, we set $\kappa = \bar{\kappa}$ and $x_{e} = x_{e}(\bar{\kappa})$ (line 32); 
\end{enumerate}

We stress that the present update rule can still be interpreted as an ``OC-like'' scheme since, for the simple case of a single constraint we are concerned with, the primal map \eqref{eq:projectionFmma} can be written out explicitly.

In our testing we have generally observed good behavior of this update rule. However, we acknowledge that the rather heuristic update for the special cases 2. and 3. may lead to oscillations in the convergence history, or even breakdowns. Since a robust optimizer is beyond the reach of a 35-lines optimization routine and beyond the purpose of this paper, the users who might face bad behavior of the \texttt{ocUpdate} are recommended to replace it with the MMA. As an example, considering the classic Matlab version of the MMA ``\texttt{mmasub}'' \citep{svanberg-note}, the buckling maximization problem, with compliance and volume constraint can be solved by replacing lines 224-226 with the following
\begin{lstlisting}[basicstyle=\scriptsize\ttfamily,breaklines=true,frame=single]
if loop==1, xOld=x(act); xOld1=xOld; low=[]; upp=[]; end
[xMMA,yy,zmma,lmid,xi,eeta,muu,zet,s1,low,upp]=mmasub_new(1,length(act),loop,x(act),max(x(act)-ocPar(1),0),...
min(x(act)+ocPar(1),1),xOld,xOld1,g0,dg0(act),0.,g1,dg1(act)',0,low,upp,a0MMA,aMMA,ccMMA,ddMMA,beta,restartAs)
xOld1 = xOld; xOld = x(act); x(act) = xMMA;
\end{lstlisting}
if the user still wants to consider a single, KS aggregated constraint. Alternatively, the following call
\begin{lstlisting}[basicstyle=\scriptsize\ttfamily,breaklines=true,frame=single]
if loop==1, xOld=x(act); xOld1=xOld; low=[]; upp=[]; end
[xMMA,yy,zmma,lmid,xi,eeta,muu,zet,s1,low,upp]=mmasub_new(2,length(act),loop,x(act),max(x(act)-ocPar(1),0),...
min(x(act)+ocPar(1),1),xOld,xOld1,g0,dg0(act),0.,g1Vec,[dg1c(act)';dg1V(act)'],0,low,upp,a0MMA,aMMA,ccMMA,...
ddMMA,beta,restartAs)
xOld1 = xOld; xOld = x(act); x( act ) = xMMA;
\end{lstlisting}
can be used for treating the two constraints separately. The parameters \texttt{a0MMA,aMMA,ccMMA} and \texttt{ddMMA} can be selected as recommended in \citep{svanberg-note}.

\begin{lstlisting}[basicstyle=\scriptsize\ttfamily,breaklines=true,frame=single]
function [x,as,lmid]=ocUpdate(loop,xT,dg0,g1,dg1,ocPar,xOld,xOld1,as,beta)
% -------------------------------- definition of asymptotes and move limits
[xU,xL] = deal(min(xT+ocPar(1),1), max(xT-ocPar(1),0));
if (loop<2.5 || restartAsy==1)
    as = xT+[-0.5,0.5].*(xU-xL)./(beta+1);
else
    tmp = (xT-xOld).*(xOld-xOld1);
    gm = ones(length(xT),1);
    [gm(tmp>0), gm(tmp<0)] = deal(ocPar(3),ocPar(2));
    as = xT + gm .* [-(xOld-as(:,1)),(as(:,2)-xOld)];
end
xL = max(0.9*as(:,1)+0.1*xT,xL);                     % adaptive lower bound
xU = min(0.9*as(:,2)+0.1*xT,xU);                     % adaptive upper bound
% ----- split (+) and (-) parts of the objective and constraint derivatives
p0_0 = (dg0>0).*dg0; q0_0 = (dg0<0).*dg0;
p1_0 = (dg1>0).*dg1; q1_0 = (dg1<0).*dg1;
[p0,q0] = deal(p0_0.*(as(:,2)-xT).^2,-q0_0.*(xT-as(:,1)).^2);
[p1,q1] = deal(p1_0.*(as(:,2)-xT).^2,-q1_0.*(xT-as(:,1)).^2);
% ---------------------- define the primal projection map and dual function
primalProj = @(lm) min(xU,max(xL,(sqrt(p0+lm*p1).*as(:,1)+sqrt(q0+lm*q1).*as(:,2))...
    ./(sqrt(p0+lm*p1)+sqrt(q0+lm*q1))));
psiDual = @(lm) g1 - ( (as(:,2)-xT)'*p1_0 - (xT-as(:,1))'*q1_0 ) + ...
    sum(p1./(max(as(:,2)-primalProj(lm),1e-12)) + q1./(max(primalProj(lm)-as(:,1),1e-12)));
% ----------------------- compute the Lagrange multiplier through bisection
lmUp = 1e6; x = xT; lmid = -1;
if psiDual( 0 ) * psiDual( lmUp ) < 0  % check if LM is within the interval
    lmid = fzero( psiDual, [ 0, lmUp ] );
    x = primalProj( lmid );                       % update desing variables
elseif psiDual(0) < 0                         % constraint cannot be active
   lmid=0; x=primalProj(lmid);
elseif psiDual( lmUp ) > 0                 % constraint cannot be fulfilled
   lmid=lmUp; x=primalProj(lmid);
end
end
\end{lstlisting}

\subsection{Some properties of the KS function}

Given the set of functions $\{g_{i}(\mathbf{x})\}_{i=1\ldots,q}$, where each $g_{i}\colon \mathbb{R}^{m}\rightarrow\mathbb{R}$ is not necessarily smooth, we can build the smooth Kreisselmeier-Steinhauser (KS) aggregation function \citep{kreisselmeier-steinhauser_79a}
\begin{equation}
 \label{eq:App-ksAggregation}
   J^{KS}[g_{i}](\mathbf{x}) = g_{\ast}(\mathbf{x}) + \frac{1}{\rho} \ln\left( \sum^{q}_{i=1}
   e^{\rho\left( g_{i}(\mathbf{x}) - g_{\ast}(\mathbf{x}) \right)} \right)
\end{equation}
depending on the parameter $\rho \in [1, \infty)$ and where $g_{\ast}(\mathbf{x}) = \max_{i=1,\ldots,q}\{g_{i}(\mathbf{x})\}$. \autoref{eq:App-ksAggregation} is a convex function if and only if all its arguments $g_{i}$, $i=1,\ldots,q$ are convex and fulfills $g_{\ast}(\mathbf{x}) \leq J^{KS}[g_{i}](\mathbf{x}) \leq g_{\ast}(\mathbf{x}) + \ln( \rho^{-1} q_{a})$, where the right inequality is tight at points where $q_{a}$ functions simultaneously attain the maximum value \citep{wrenn_89a,raspanti-etal_00a}.

With the aggregation of several constraint functions in mind, the following considerations are in order. At a given point $\mathbf{x}$, if all $g_{i}(\mathbf{x}) > 0$, then $g_{\ast}(\mathbf{x}) > 0$ and $(g_{i} - g_{\ast})(\mathbf{x}) \leq 0$ for all $i$. If all $g_{i}(\mathbf{x}) < 0$, then $g_{\ast}(\mathbf{x}) < 0$ and $(g_{i} - g_{\ast})(\mathbf{x}) \leq 0$ for all $i$. In both cases $J^{KS}[g_{i}](\mathbf{x}) \geq g_{\ast}(\mathbf{x})$ and, if all the terms are negative this means that $J_{KS}[g_{i}](\mathbf{x})$ is closer to zero than all $g_{i}(\mathbf{x})$. Thus it gives an upper bound to the constraint that is closer to become active. The interesting case is when some $g_{i}(\mathbf{x}) > 0$ and some other $g_{j}(\mathbf{x}) < 0$. Obviously, $g_{\ast}(\mathbf{x}) > 0$ and $J^{KS}[g_{i}](\mathbf{x})$ gives an upper bound to the positive maximum; in other words, to the most violated constraint.

\section{Matlab code}
 \label{App:matlabCodeBuckling}

 \begin{lstlisting}
function topBuck250(nelx,nely,penalK,rmin,ft,ftBC,eta,beta,maxit,ocPar,Lx,penalG,nEig,pAgg,prSel,x0)
% ---------------------------- PRE. 1) MATERIAL AND CONTINUATION PARAMETERS
[E0,Emin,nu] = deal(1,1e-6,0.3);                                           % Young's moduli & Poisson's ratio
penalCntK = {25,1,25,0.25};                                                % continuation scheme on K-penal
penalCntG = {25,1,25,0.25};                                                % " " on G-penal
betaCnt  = { 400,24,25,2};                                                 % " " on beta
pAggCnt  = { 2e5,1,25,2};                                                  % " " on the KS aggregation factor
cnt = @(v,vCn,l) v+(l>=vCn{1}).*(v<vCn{2}).*(mod(l,vCn{3})==0).*vCn{4};    % function applying continuation
if prSel{1}(1) == 'V', volfrac = 1.0; else, volfrac = prSel{2}(end); end   % initialize volume fraction
% ----------------------------------------- PRE. 2) DISCRETIZATION FEATURES
Ly = nely/nelx*Lx;                                                         % recover Ly from aspect ratio
nEl = nelx*nely;                                                           % number of elements
elNrs = reshape(1:nEl,nely,nelx);                                          % element numbering
nodeNrs = int32(reshape(1:(1+nely)*(1+nelx),1+nely,1+nelx));               % node numbering (int32)
cMat = reshape(2*nodeNrs(1:end-1,1:end-1)+1,nEl,1)+int32([0,1,2*nely+[2,3,0,1],-2,-1]);% connectivity matrix
nDof = (1+nely)*(1+nelx)*2;                                                % total number of DOFs
% ---------------------------------------------- elemental stiffness matrix
c1 = [12;3;-6;-3;-6;-3;0;3;12;3;0;-3;-6;-3;-6;12;-3;0;-3;-6;3;12;3;...
    -6;3;-6;12;3;-6;-3;12;3;0;12;-3;12];
c2 = [-4;3;-2;9;2;-3;4;-9;-4;-9;4;-3;2;9;-2;-4;-3;4;9;2;3;-4;-9;-2;...
    3;2;-4;3;-2;9;-4;-9;4;-4;-3;-4];
Ke = 1/(1-nu^2)/24*(c1+nu.*c2);                                            % lower symmetric part of Ke
Ke0(tril(ones(8))==1) = Ke';
Ke0 = reshape(Ke0,8,8);
Ke0 = Ke0+Ke0'-diag(diag(Ke0));                                            % recover full elemental matrix
[sI,sII] = deal([]);
for j = 1:8      % build assembly indices for the lower symmetric part of K
    sI = cat(2,sI,j:8);
    sII = cat(2,sII, repmat(j,1,8-j+1));
end
[iK,jK] = deal(cMat(:,sI)',cMat(:,sII)');
Iar = sort([iK(:),jK(:)],2,'descend');                                     % indices for K assembly
if any(prSel{1}=='B') % >>>>>>>>>>>>>>>> PERFORM ONLY IF BUCKLING IS ACTIVE #B#
    Cmat0 = [1,nu,0;nu,1,0;0,0,(1-nu)/2]/(1-nu^2);                         % non-dimensional elasticity matrix
    xiG = sqrt(1/3)*[-1,1]; etaG = xiG; wxi = [1,1]; weta = wxi;           % Gauss nodes and weights
    xe = [-1,-1;1,-1;1,1;-1,1].*Lx/nelx/2;                                 % dimensions of the elements
    lMat = zeros(3, 4); lMat(1, 1) = 1; lMat(2, 4) = 1; lMat(3, 2:3) = 1;  % placement matrix
    dN = @(xi,zi) 0.25*[zi-1,1-zi,1+zi,-1-zi; xi-1,-1-xi,1+xi,1-xi];       % shape funct. logical derivatives
    B0 = @(gradN) lMat * kron(gradN,eye(2));                               % strain-displacement matrix
    [indM,t2ind] = deal([1,3,5,7,16,18,20,27,29,34],[ 2,3,4,6,7,9 ]);      % auxiliary set of indices (1)
    [iG,jG] = deal(iK(indM,:),jK(indM,:));                                 % indexing of unique G coefficients
    IkG = sort([iG(:), jG(:)],2,'descend');                                % indexing G entries (lower half)
    [a1,a2]=deal(reshape(IkG(:,2),10,nEl)', reshape(IkG(:,1),10,nEl)');    % auxiliary set of indices (2)
    dZdu = zeros(10,8);                                                    % build U-derivative of matrix G
    for ii = 1 : 8                    % loop on the displacement components
        tt = 0; Uvec = zeros(8,1); Uvec(ii,1) = 1;                         % set a single displ. component
        se = Cmat0*B0((dN(0,0)*xe)\dN(0,0))*Uvec;                          % stresses at the element center
        for j = 1 : length(xiG)
            for k = 1 : length(etaG)
                xi = xiG(j); zi = etaG(k);                                 % current integration points
                w = wxi(j)*weta(k)*det(dN(xi,zi)*xe);                      % current integration weight
                gradN = (dN(xi,zi)*xe)\dN(xi,zi);                          % shape funct. physical derivatives
                B1 = [kron(gradN,[1,0]); kron(gradN,[0,1])];               % deformation gradient
                tt = tt+(B1'*kron(eye(2),[se(1),se(3);se(3),se(2)])*B1)*w; % current contribution to dG/du_i
            end
        end
        dZdu(:,ii) = tt([1,3,5,7,19,21,23,37,39,55])';                     % extract independent coefficients
    end
    dZdu(t2ind,:) = 2*dZdu(t2ind,:);                                       % x2 columns for v-m-v product
    fKS = @(p,v)max(v)+log(sum(exp(p*(v-max(v)))))/p;                      % KS aggregation function
    dKS = @(p,v,dv)sum(exp(p.*(v-max(v)))'.*dv,2)./sum(exp(p.*(v-max(v))));% derivative of the KS function
end % <<<<<<<<<<<<<<<<<<<<<<<<<<<<<<<<<<<<<<<<<<<<<<<<<<<<<<<<<<<<<<<<<<<<< #B#
% ----------------------------- PRE. 3) LOADS, SUPPORTS AND PASSIVE DOMAINS
fixed = 1:2*(nely+1);                                                      % restrained DOFs (cantilever)
lcDof = 2*nodeNrs(nely/2+1+[-8:8],end)-1;                                  % loaded DOFs
modF = 1e-3/Ly/(length(lcDof)-1);                                          % modulus of the force density
F = fsparse(lcDof,1,-modF,[nDof,1]);                                       % define load vector
[F(lcDof(1)),F(lcDof(end))] = deal(F(lcDof(1))/2,F(lcDof(end))/2);         % consistent load on end nodes
[pasS,pasV] = deal(elNrs(nely/2+[-9:10],end-9:end),[]);                    % define passive domains
free = setdiff(1:nDof, fixed);                                             % set of free DOFs
act = setdiff((1:nEl)',union(pasS(:),pasV(:)));                            % set of active design variables
% ------------------------- PRE. 4) PREPARE FILTER AND PROJECTION OPERATORS
if ftBC == 'N', bcF = 'symmetric'; else, bcF = 0; end                      % select filter BC
[dy,dx] = meshgrid(-ceil(rmin)+1:ceil(rmin)-1,-ceil(rmin)+1:ceil(rmin)-1);
h = max(0,rmin-sqrt(dx.^2+dy.^2));                                         % convolution kernel
Hs = imfilter(ones(nely,nelx),h,bcF);                                      % matrix of weights
dHs = Hs;
prj = @(v,eta,beta) (tanh(beta*eta)+tanh(beta*(v(:)-eta)))./...
    (tanh(beta*eta)+tanh(beta*(1-eta)));                                   % relaxed Heaviside projection
deta = @(v,eta,beta) -beta*csch(beta).*sech(beta*(v(:)-eta)).^2 .* ...
    sinh(v(:)*beta).*sinh((1-v(:))*beta);                                  % projection eta-derivative
dprj = @(v,eta,beta) beta*(1-tanh(beta*(v-eta)).^2)./...
    (tanh(beta*eta)+tanh(beta*(1-eta)));                                   % projection x-derivative
% ------------------------ PRE. 5) ALLOCATE AND INITIALIZE OTHER PARAMETERS
[x,dsK,dsG,dmKS,dV] = deal(zeros(nEl,1));                                  % initialize vectors of size nElx1
[phiDKphi,phiDGphi,adj] = deal(zeros(nEl,nEig));                           % " " of size nElxnEig
U = zeros(nDof,1); phi = zeros(nDof,nEig); adjL = phi; adjV = phi;         % " " of size nDofx1 & nDofxnEig
dV(act,1) = 1/nEl;                                                         % derivative of volume fraction
[xpOld,loop,restartAs,ch,plotL,plotR,muVec] = deal(0,0,0,1,[],[],[]);      % misc array & parameters
if nargin > 15
    load(x0); x = xInitial;                                                % initialize design from saved data
else
    x(act) = (volfrac*(nEl-length(pasV))-length(pasS))/length(act);        % volume fraction on "active" set
    x(pasS) = 1;                                                           % set x=1 on "passive solid" set
end
xPhys = x; clear iK jK iG jG dx dy;                                        % initialize xPhys and free memory
%% ________________________________________________ START OPTIMIZATION LOOP
while loop < maxit && ch > 1e-6
  loop = loop+1;                                                           % update iteration counter
  % --------------------------------- RL. 1) COMPUTE PHYSICAL DENSITY FIELD
  xTilde = imfilter(reshape(x,nely,nelx),h,bcF)./Hs;                       % compute filtered field
  xPhys(act) = xTilde(act);                                                % modify active elements only
  if ft > 1                                                                % apply projection
      f = (mean(prj(xPhys,eta,beta))-volfrac)*(ft==3);                     % function (volume of x-projected)
      while abs(f) > 1e-6 && prSel{1}(1) ~= 'V'                            % Newton loop for finding opt. eta
          eta = eta-f/mean(deta(xPhys(:),eta,beta));
          f = mean(prj(xPhys,eta,beta))-volfrac;
      end
      dHs = Hs./reshape(dprj(xPhys,eta,beta),nely,nelx);                   % modification of the sensitivity
      xPhys = prj(xPhys,eta,beta);                                         % compute projected field
  end
  ch = max(abs(xPhys-xpOld)); xpOld = xPhys;
  % -------------------------- RL. 2) SETUP AND SOLVE EQUILIBRIUM EQUATIONS
  sK = (Emin+xPhys.^penalK*(E0-Emin));                                     % stiffness interpolation
  dsK(act) = penalK*(E0-Emin)*xPhys(act).^(penalK-1);                      % derivative of " "
  sK = reshape(Ke(:)*sK',length(Ke)*nEl,1);
  K = fsparse(Iar(:,1),Iar(:,2),sK,[nDof,nDof]);                           % assemble stiffness matrix
  K = K+K'-diag(diag(K));                                                  % symmetrization of K
  dK = decomposition(K(free,free),'chol','lower');                         % decompose K and store factor
  U(free) = dK \ F(free);                                                  % solve equilibrium system
  dc = -dsK.*sum((U(cMat)*Ke0).*U(cMat),2);                                % compute compliance sensitivity
  if any(prSel{1}=='B') % >>>>>>>>>>>>>> PERFORM ONLY IF BUCKLING IS ACTIVE #B#
  % ---------------------------------- RL. 3) BUILD STRESS STIFFNESS MATRIX 
  sGP = (Cmat0*B0((dN(0,0)*xe)\dN(0,0))*U(cMat)')';                        % stresses at elements centroids
  Z = zeros(nEl,10);      % allocate array for compact storage of Ge coeff.
  for j = 1:length(xiG)                       % loop over quadrature points
    for k = 1:length(etaG)
        % ---------------------------- current integration point and weight
        xi = xiG(j); zi = etaG(k); w = wxi(j)*weta(k)*det(dN(xi,zi)*xe);
        % - reduced represenation of strain-displacement matrix (see paper)
        gradN = (dN(xi,zi)*xe)\dN(xi,zi);                                  % shape funct. physical derivatives
        a = gradN(1,:); b = gradN(2,:); B = zeros(3,10);
        l = [1,1;2,1;3,1;4,1;2,2;3,2;4,2;3,3;4,3;4,4];
        for jj = 1:10
            B(:,jj) = [a(l(jj,1))*a(l(jj,2)); ...
                       b(l(jj,1))*b(l(jj,2)); ...
                       b(l(jj,2))*a(l(jj,1))+b(l(jj,1))*a(l(jj,2))]; 
        end
        % ----------- current contribution to (unique ~= 0) elements of keG
        Z = Z+sGP*B*w;
    end
  end
  sG0 = E0*xPhys.^penalG;                                                  % stress interpolation
  dsG(act) = penalG*E0*xPhys(act).^(penalG-1);                             % derivative of " "
  sG = reshape((sG0.*Z)',10*nEl,1);
  G = fsparse(IkG(:,1)+1,IkG(:,2)+1,sG,[nDof,nDof])+...
         fsparse(IkG(:,1),  IkG(:,2),  sG,[nDof,nDof]);                    % assemble global G matrix
  G = G+G'-diag(diag(G));                                                  % symmetrization of G
  % ------------------------------ RL. 4) SOLVE BUCKLING EIGENVALUE PROBLEM
  matFun = @(x) dK\(G(free,free)*x);                                       % matrix action function
  [eivecs,D] = eigs(matFun,length(free),nEig+4,'sa');                      % compute eigenvalues
  [mu,ii] = sort(diag(-D),'descend');                                      % sorting of eigenvalues (mu=-D(i))
  eivSort = eivecs(:,ii(1:nEig));                                          % sort eigenvectors accordingly
  phi(free,:) = eivSort./sqrt(diag(eivSort'*K(free,free)*eivSort)');       % orthonormalize (phi'*K*phi=1)
  % ----------------------------------- RL. 5) SENSITIVITY ANALYSIS OF BLFs
  dkeG = dsG.*Z;                                                           % x-derivative of Ge
  dkeG(:,t2ind) = 2*dkeG(:,t2ind);                                         % x2 columns for v-m-v product
  for j = 1:nEig     % loop on the eigenvalues included in the optimization
      % 1) ------ Term due to the elastic stiffness matrix (phi'*dK/dx*phi)
      t = phi(:,j);
      phiDKphi(:,j) = dsK.*sum((t(cMat)*Ke0).*t(cMat),2);
      % 2) -------------- Term due to the geometric matrix (phi'*dG/dx*phi)
      p = t(a1).*t(a2)+t(a1+1).*t(a2+1);
      phiDGphi(:,j) = sum(dkeG.*p,2);
      % 3) ----------------------------------------- Setup of adjoint loads
      tmp = zeros(nDof,1);
      for k = 1:8              % contribution of each term dKg/du_i, i=1:nD
          tmp(cMat(:,k)) = tmp(cMat(:,k))+(sG0.*p)*dZdu(:,k);
      end
      adjL(:,j) = tmp;
  end
  % ----------- solve the adjoint problem and compute the term (U'*dK/dx*V)
  adjV(free,:) = dK \ adjL(free,:);               % use the stored K factor
  for j = 1 : nEig
      vv = adjV(:,j);
      adj(:,j) = dsK.*sum((U(cMat)*Ke0).*vv(cMat),2);
  end 
  % --------------- overall sensitivity expression for the "mu" eigenvalues
  dmu = -(phiDGphi+mu(1:nEig )'.*phiDKphi-adj);
  end % <<<<<<<<<<<<<<<<<<<<<<<<<<<<<<<<<<<<<<<<<<<<<<<<<<<<<<<<<<<<<<<<<<<
  % ---------------------- RL. 6) SELECT OBJECTIVE FUNCTION AND CONSTRAINTS
  if loop==1, c0=F'*U; v0=mean(xPhys(:)); end % initial compliance & volume fraction
  switch prSel{1}                    % select optimization problem to solve
      case ['C','V']           % minimize compliance with volume constraint
          g0 = F'*U/c0;
          dg0 = imfilter(reshape(dc/c0,nely,nelx)./dHs,h,bcF);
          g1 = mean(xPhys(:))/volfrac-1;
          dg1 = imfilter(reshape(dV/volfrac,nely,nelx)./dHs,h,bcF);
      case ['V','C']           % minimize volume with compliance constraint
          g0 = mean(xPhys(:))./v0;
          dg0 = imfilter(reshape(dV/v0,nely,nelx)./dHs,h,bcF);
          g1 = (F'*U)/(prSel{2}*c0)-1;
          dg1 = imfilter(reshape(dc/(prSel{2}*c0),nely,nelx)./dHs,h,bcF);
      case ['B','C','V']% maximize BLF with compliance & volume constraints (Eq. 13 paper)
          if loop==1, muKS0=fKS(pAgg,mu(1:nEig)); g0=1; cMax=prSel{2}(1);
          else, g0=fKS(pAgg,mu(1:nEig))/muKS0; end                         % KS aggregation of mu (=1/lambda)
          dmKS = dKS(pAgg,mu(1:nEig),dmu);                                 % KS aggregation of dmu
          dg0 = imfilter(reshape(dmKS/muKS0,nely,nelx)./dHs,h,bcF);        % back-filter KS sensitivity
          % -- Constraint function: KS aggregation of compliance and volume
          g1Vec = [F'*U;mean(xPhys(:))]./[cMax*c0;volfrac]-1;              % set of constraints ['C','V']
          dg1c = imfilter(reshape(dc/(cMax*c0),nely,nelx)./dHs,h,bcF);     % back-filter compliance derivative
          dg1V = imfilter(reshape(dV/volfrac,nely,nelx)./dHs,h,bcF);       % back-filter volume derivative
          g1 = fKS(pAgg,g1Vec);                                            % aggregate the two constraints
          dg1 = dKS(pAgg,g1Vec,[dg1c(:),dg1V(:)]);                         % sensitivity of the KS constraint
          plotL(loop,:) = [1/g0/muKS0,1/mu(1)]; strL='KS(-),\lambda_1(--)';
          plotR(loop,:) = [g1,g1Vec']; strR='g_1(-),gC(--),gV(.-)';
          muVec(loop,:) = mu';
      case ['V','C','B']    % min volume with compliance & BLF constraints (Eq. 14 paper)
          g0 = mean(xPhys(:))./v0;
          dg0 = imfilter(reshape(dV/volfrac,nely,nelx)./dHs,h,bcF);
          % ---- Constraint function: KS aggregation of BLFs and compliance
          muKS = fKS(pAgg,mu(1:nEig));                                     % KS aggregation of mu
          dmKS = dKS(pAgg,mu(1:nEig),dmu);                                 % KS aggregation of dmu
          g1Vec = [prSel{2}(2)*muKS;F'*U]./[1;prSel{2}(1)*c0]-1;           % set of constraints 'B','C'
          dg1l = imfilter(reshape(dmKS*prSel{2}(2),nely,nelx)./dHs,h,bcF); % back-filter dmu
          dg1c = imfilter(reshape(dc/(prSel{2}(1)*c0),nely,nelx)./dHs,h,bcF);% back-filter dc
          g1 = fKS(pAgg,g1Vec);                                            % aggregate the two constraints
          dg1 = dKS(pAgg,g1Vec,[dg1l(:),dg1c(:)]);                         % sensitivity of the KS constraint
          plotL(loop,:) = g0; strL = 'g_0';
          plotR(loop,:) = [g1,g1Vec']; strR='g_1(-),gL(--),gC(.-)';
          muVec = cat(1,muVec,mu');
  end
  % ---------------------------------------- RL. 7) UPDATE DESIGN VARIABLES
  if loop==1, xOld = x(act); xOld1 = xOld; as = []; end                    % initialize MMA history parameters
  [x0,as,lmid]=ocUpdate(loop,x(act),dg0(act),g1,dg1(act),ocPar,xOld,xOld1,as,beta,restartAs);
  xOld1 = xOld; xOld = x(act); x(act) = x0;
  % ----------------------------------------- RL. 8) PRINT AND PLOT RESULTS 
  fprintf('It.:%2i g0:%7.4f g1:%0.2e penalK:%7.2f penalG:%7.2f eta:%7.2f beta:%7.1f ch:%0.3e lm:%0.3e\n', ...
    loop,g0,g1,penalK,penalG,eta,beta,ch,lmid);
  if any(prSel{1} == 'B')  % plot design, g0 & g1 evolution, BLFs evolution
      subplot(2,2,1:2);
      colormap(gray); imagesc(1-reshape(xPhys,nely,nelx));
      caxis([0,1]); axis equal; axis off; drawnow; title('Current design');
      subplot(2,2,3)
      yyaxis left; plot(1:loop,plotL); ylabel(strL);
      yyaxis right; plot(1:loop,plotR); ylabel(strR); title('Objective and constraint');
      subplot(2,2,4)
      plot(1:loop,1./muVec(:,1:4)); title('Lowest BLFs');
  else                                       % plot the current design only
      colormap(gray); imagesc(1-reshape(xPhys,nely,nelx));
      caxis([ 0,1]); axis equal; axis off; drawnow;
  end
  %  apply continuation on penalization(s), beta & aggregation parameter(s)
  penalKold = penalK; penalGold = penalG; betaOld = beta;
  [penalK,penalG,beta,pAgg] = deal(cnt(penalK, penalCntK, loop), ...
      cnt(penalG,penalCntG,loop), cnt(beta,betaCnt,loop), cnt(pAgg,pAggCnt,loop));
  if (beta-betaOld~= 0 || penalK-penalKold~=0 || penalG-penalGold~=0)
      restartAs = 1; else, restartAs = 0; end                              % restart asymptotes if needed
end
end
\end{lstlisting}


\clearpage

\nomenclature{$m$}{number of elements (\texttt{nEl})}
\nomenclature{$n$}{number of elements (\texttt{nDof})}
\nomenclature{$\mathbf{x}$}{vector of design variables (\texttt{x})}
\nomenclature{$\tilde{\mathbf{x}}$}{vector of filtered variables (\texttt{xTilde})}
\nomenclature{$\hat{\mathbf{x}}$}{vector of physical variables (\texttt{xPhys})}

\nomenclature{$\eta$}{projection threshhold (\texttt{eta})}
\nomenclature{$\beta$}{projection sharpness factor (\texttt{beta})}
\nomenclature{$r_{\rm min}$}{density filter radius (\texttt{rmin})}
\nomenclature{$E_{1}$}{(\texttt{E0})}
\nomenclature{$E_{0}$}{(\texttt{Emin})}
\nomenclature{$p_{K}$}{(\texttt{penalK})}
\nomenclature{$p_{G}$}{(\texttt{penalG})}
\nomenclature{$f$}{(\texttt{volfrac})}
\nomenclature{$c$}{(\texttt{c})}
\nomenclature{$\mathbf{F}$}{(\texttt{F})}
\nomenclature{$\mathbf{u}$}{(\texttt{U})}
\nomenclature{$K$}{(\texttt{K})}
\nomenclature{$G$}{(\texttt{G})}
\nomenclature{$\mu$}{(\texttt{mu})}
\nomenclature{$\boldsymbol{\varphi}$}{(\texttt{phi})}

\printnomenclature
\end{document}